\documentclass[12pt,a4paper]{article}

\usepackage[latin1]{inputenc}
\usepackage{amsmath, amsfonts, amssymb, graphicx, amsthm}
\usepackage{graphicx}
\usepackage[spaces,hyphens]{url}
\makeindex
\newtheorem{theorem}{Theorem}[section]
\newtheorem{corollary}[theorem]{Corollary}
\newtheorem{defi}[theorem]{Definition}
\newtheorem{example}[theorem]{Example}
\newtheorem{lemma}[theorem]{Lemma}

\newtheorem{proposition}[theorem]{Proposition}
\newtheorem{remark}[theorem]{Remark}

\newcommand*{\N}{\mathbb{N}}
\newcommand*{\R}{\mathbb{R}}
\newcommand*{\C}{\mathbb{C}}
\newcommand*{\Bcal}{\mathcal{B}}

\newcommand*{\tr}{\text{{\normalfont{tr}}}\,}

\newcommand*{\Id}{\text{\normalfont{Id}}}

\newcommand{\braket}[2]{\langle\,#1\,\vert\,#2\,\rangle}
\newcommand{\bra}[1]{\langle\,#1\,\rvert}
\newcommand{\ket}[1]{\lvert\,#1\,\rangle}

\newcommand{\norm}[1]{\left\Arrowvert#1\right\Arrowvert}
\newcommand{\abs}[1]{\left\arrowvert#1\right\arrowvert}

\newcommand{\inner}[2]{\langle\,#1\,,\,#2\,\rangle}

\newcommand{\matrixx}[2][cc]{\left(\begin{array}{#1}#2\end{array}\right)}

\newcommand{\dm}{\,d\mu}
\newcommand{\dn}{d\nu}
\renewcommand{\Pr}{\mathbb{P}}
\newcommand*{\E}{\mathbb{E}}
\newcommand*{\supp}{\text{supp}\,}

\title{Thermodynamic Formalism for Quantum Channels: Entropy, Pressure, Gibbs Channels and generic properties}

\author{Jader E. Brasil, Josu\'e Knorst and Artur O. Lopes}

\begin{document}
\maketitle
\centerline{Abstract}
\smallskip

Denote  $M_k$ the set of complex $k$ by $k$ matrices.
We will analyze here quantum channels $\phi_L$ of the following kind: given a measurable function  $L:M_k\to M_k$ and the measure $\mu$ on $M_k$ we define the linear operator $\phi_L:M_k \to M_k$, via the expression $\rho \,\to\,\phi_L(\rho) = \int_{M_k} L(v) \rho {L(v)}^\dagger \, \dm(v)$.

A recent paper by T. Benoist, M. Fraas, Y. Pautrat, and C. Pellegrini is our starting point. They considered the case where $L$ was the identity.

Under some mild assumptions on the quantum channel  $\phi_L$
we analyze
the eigenvalue property for $\phi_L$ and we  define entropy for such channel.
For a fixed $\mu$ (the \textit{a priori} measure) and for a given a Hamiltonian $H: M_k \to M_k$ we present a version of the Ruelle Theorem: a variational principle  of pressure (associated to such $H$) related to an eigenvalue problem for the Ruelle operator. We introduce the concept of Gibbs channel.

We also show that for a fixed $\mu$ (with more than one point in the support) the set of $L$ such that it is $\phi$-Erg (also irreducible) for $\mu$ is a generic set.

We describe a related process $X_n$, $n\in \mathbb{N}$, taking values on the projective space $ P(\C^k)$ and analyze the question of the existence of invariant probabilities.
We also consider an associated process $\rho_n$, $n\in \mathbb{N}$,   with values on $\mathcal{D}_k$  ($\mathcal{D}_k$
is the set of density operators). Via the  barycenter, we associate the invariant probability  mentioned above with  the  density operator fixed for $\phi_L$.

\smallskip

\section{Introduction}

There are many different definitions  and meanings for the concept of quantum dynamical entropy. We mention first the more well-known concepts due to  Connes-Narnhofer-Thirring (see $\cite{CNT}$), Alicki-Fannes (see $\cite{AF}$), Accardi-Ohya-Watanabe (see $\cite{Acar}$), Stormer \cite{Stor}
and Kossakowski-Ohya-Watanabe (see $\cite{KOW}$). In this case, the entropy  can be exactly computed for several examples of quantum
dynamical systems.

A different approach appears in $\cite{SSZ}$ and $\cite{SZ}$ where the authors present their definition of
quantum dynamical entropy (see also $\cite{AW}$).

Classical texts on quantum entropy are $\cite{AF1}$, $\cite{Ben}$, $\cite{Ben1}$ and $\cite{OP}$, and for 
quantum channels we also mention $\cite{JP}$, $\cite{Li}$,  and $\cite{wolf2012quantum}$.

We present here a certain concept of dynamical quantum entropy. A confirmation that this entropy is
in fact a concept that describes valuable information from a dynamic point of view is its relationship with Lyapunov exponents as presented in \cite{BKL} by the same authors. Lyapunov exponents are quite important tools that are used  in Physics, Dynamics, and Fractals. Moreover,  in  \cite{BKL} we will  show that the purification property is  $C^0$-generic.

One of the most challenging open problems in quantum information theory is to introduce a good definition capable of quantifying how \textit{entanglement} behaves when part of an entangled state is sent through a quantum channel. Therefore the understanding of quantum channels is a problem of central importance.

Denote  $M_k$ the set of complex $k$ by $k$ matrices.
We will analyze here quantum channels $\phi_L$ of the following kind: given a measurable function  $L:M_k\to M_k$ and the measure $\mu$ on $M_k$ we define the linear operator $\phi_L:M_k \to M_k$, via the expression $\rho \,\to\,\phi_L(\rho) = \int_{M_k} L(v) \rho {L(v)}^\dagger \, \dm(v)$.

The probability $\mu$ will play the role of an {\it a priori} probability for defining entropy (in the spirit of \cite{LMMS}) as described in section \ref{ent}.

In $\cite{benoist2017invariant}$ the authors present interesting results for the case $L=I$.
This paper is our starting point and we follow its notation as much as possible. Given $L$ (as above) one can consider in the setting of $\cite{benoist2017invariant}$  a new probability $\mu_L= \mu \circ L^{-1}$ and part of the results presented here can be recovered from there (using $\mu_L$ instead of $\mu$).

We will present all the proofs here using $L$ and  $\mu$ as above (and not via $\mu_L$) because this will be more natural for our future reasoning (for instance when analyzing generic properties).

 In the Thermodynamic Formalism version of Quantum Information, the $L$ will help on the one hand to express the {\it analogous concept of function (even the analog of a Hamiltonian)} and on the other hand, a certain class of $L$ - together with the {\it a priori} probability $\mu$ on $M_k$ - will  help to describe the analogous concept of {\it invariant probability}. Later we will elaborate on that.

This paper is self-contained.

For a fixed $\mu$ and a general $L$ we present a natural concept of entropy for a channel in order to develop a version of  Gibbs formalism
which seems natural to  us.  Example  $\ref{Mark1}$ in Section $\ref{exam}$  (the Markov model in quantum information) will show that our definition is a natural extension of the classical  concept of entropy.  We point out that
the definition  of entropy we will consider here is a generalization of the concept described on the papers $\cite{baraviera2010thermodynamic}$, $\cite{BLLT2}$ and $\cite{BLLT1}$.
This particular way of defining entropy is inspired by the results of
$\cite{slomczynski2003dynamical}$ which consider iterated function systems.

For  a given  $H: M_k \to M_k$ (which plays the role of a Hamiltonian)  we present a version of the Ruelle Theorem for $\phi_H$: a variational principle  of pressure related to an eigenvalue problem for a kind of Ruelle operator (see Theorem \ref{peqlog}).

A question of terminology:  the operator $H$ (mentioned above as Hamiltonian) could also be naturally called Liouvillian; it would make perfect sense taking into account that $M_k$  is an algebra of quantum observables where the operator
acts (Heisenberg picture of QM). The notation  $L$ used by the authors  in  \cite{benoist2017invariant} was
probably inspired by their  understanding that $L$ plays  the role of a Liouvillian operator.

We say that  $E\subset \C^k$ is $(L, \mu)$-invariant if $L(v)(E)\subset E$, for all $v  $ in the  support of $\mu$.
Given $L: M_k\to M_k$ and  $\mu$ on $M_k$, we say that $L$ is $\phi$-Erg for $\mu$, if there exists an unique minimal non-trivial space $E$, such that, $E$ is $(L, \mu)$-invariant.
We will show in Section $\ref{gen}$ that for a fixed $\mu$ (with more than one point in the support) the set of $L$ such that it is $\phi$-Erg for $\mu$ is generic. In fact, the set of $L$ which are irreducible is dense according to Theorem $\ref{coi}$.

The introduction of this variable $L$ allows us to consider questions of a generic nature in this type of problem.


We point out that here we explore the point of view that the (discrete-time dynamical) classical Kolmogorov-Shannon entropy of an invariant probability is in some way
attached to an \textit{a priori} probability (even if this is not transparent on the
classical definition). This point of view becomes more clear when someone tries to
analyze the generalized $XY$ model (the symbolic space $M^\mathbb{N}$ where the alphabet  $M$ is  a compact metric space) which is a case with the property  that each point has an uncountable number of
preimages (see $\cite{LMMS}$ and $\cite{BCL}$ for discussion). In the dynamical setting of $\cite{LMMS}$
to define entropy it is necessary first to introduce  the transfer (Ruelle) operator
(which we claim - in some sense - is a more fundamental concept than entropy)
which requires an \textit{a priori} probability (not a general measure).  Our results correspond to the case where  the alphabet (that in some sense corresponds to the support of the {\it a priori} probability $\mu$)
can be uncountable.

The point of view of defining entropy via the limit of dynamical partitions is not suitable for the generalized $X\,Y$ model. We are just saying that in any case the concept of entropy can be recovered via the Ruelle operator.

We point out, as a curiosity, that for the computation of the classical
Kolmogorov-Shannon entropy of a shift invariant probability on ${\{1,2,\ldots,d\}}^\mathbb{N}$
one should  take as the \textit{a priori} measure (not a probability) the counting
measure on $\{1,2,\ldots,d\}$ (see discussion in $\cite{LMMS}$).  In the case, we take as
\textit{a priori probability} $\mu$ the uniform normalized probability on
$\{1,2,\ldots,d\}$ the entropy will be negative (it will be Kolmogorov-Shannon entropy
- $\log d$). In this case the independent $1/d$ probability on ${\{1,2,\ldots,d\}}^{\mathbb{N}}$ will have maximal entropy equal $0$.

A general reference for Thermodynamic Formalism is \cite{PP} and \cite{LNotes}.
\medskip

We point out that we consider here Quantum Channels but the associated  discrete-time process is associated with a Classical Stochastic Process (a probability on the infinite product of an uncountable state space) and not to a quantum spin-lattice, where it is required the use of the tensor product (see \cite{LMMM} and \cite{BLMM}).

After some initial sections describing basic properties which will be required later
we analyze in Section $\ref{eig}$
the eigenvalue property for $\phi_L$.

Under some mild assumptions on $\phi_L$, we define the entropy of the channel $\phi_L$ in Section $\ref{ent}$.
For a fixed $\mu$ (the \textit{a priori} measure) and a given  Hamiltonian
$H: M_k \to M_k$ we present a variational principle  of pressure and we associate with all
this an eigenvalue problem on Section $\ref{eig}$. In Definition $\ref{ort}$ we introduce the concept of Gibbs channel for the Hamiltonian $H$ (or,  for the channel $\phi_H$).

In Section  $\ref{pro1}$ we describe (adapting $\cite{benoist2017invariant}$ to the
present setting) a process $X_n$, $n\in \mathbb{N}$, taking values on the projective space $ P(\C^k)$. We also analyze the existence of an initial invariant probability for this process (see Theorem $\ref{existence_inv_prob}$).

In Section $\ref{pro2}$ we consider a process $\rho_n$, $n\in \mathbb{N}$  (called quantum trajectory by T. Benoist, M. Fraas, Y. Pautrat, and C. Pellegrini) taking values on $\mathcal{D}_k$, where $\mathcal{D}_k$
	is the set of density operators on $M_k$. Using the definition of barycenter taken from $\cite{slomczynski2003dynamical}$ we relate in proposition $\ref{barycenter}$ the invariant probabilities of Section $\ref{pro1}$ with the fixed point of Section $\ref{eig}$.

In Section $\ref{gen}$,  for a fixed measure $\mu$, we show that $\phi$-Erg (and  also irreducible) is a generic property for $L$ (see Corollary $\ref{corimp}$).

In Section $\ref{exam}$, we present several examples that will  help the reader in understanding the theory.  Example $\ref{Mark1}$ shows that the definition of entropy for Quantum Channels described here is the natural generalization of the classical concept of entropy. In another example in this section, we consider the case where $\mu$ is a probability with support on a linear space of $M_2$ (see Example $\ref{gre}$), and among other things we estimate the entropy of the channel.

In the final section \ref{Con} we will present some clarifications on which directions our work is related to relevant issues in the area connected to quantum entropy.

We thank the referees for the careful reading of our manuscript and for providing us with many suggestions to improve the reading of the paper.

\smallskip

\medskip

\section{General properties}

\hspace{13pt} We present some basic definitions.

We denote by $M_k$, $k \in \mathbb{N}$,  the set of complex $k$ by $k$ matrices. We consider $\mathcal{M}$ the standard Borel sigma-algebra over $M_k$ and on $\mathbb{C}^k$, we consider the canonical Euclidean inner product.

We denote by $\Id_k$ the identity matrix on $M_k$.

According to our notation, $\dagger$ denotes the operation  of taking the dual of a matrix with respect to the canonical inner product on $\C^k$.

Here tr denotes the trace of a matrix.

Given two matrices $A$ and $B$ we define the Hilbert-Schmidt product
\[\inner{A}{B} \,\,= \, \text{ tr}\,\,  (A\, B^\dagger).\]

This induces a  norm $\norm{A}=\sqrt{ \inner{A}{A} }$ on the Hilbert space $M_k$ which will be called the Hilbert-Schmidt norm.

Given a linear operator $\phi$ on $M_k$ we denote by $\phi^*: M_k \to M_k$ the dual linear operator in the sense of Hilbert-Schmidt, that is, if for all $X,Y$ we get

\[\inner{ \phi(X)}{Y}\, = \inner{X}{\phi^* (Y)}.\]

Now, consider a measure $\mu$ on $\mathcal{M}$.

For an integrable transformation  $F: M_k \to M_k$:
\[\int_{M_k} F(v) \, \dm(v) = {\left( \int_{M_k} {F(v)}_{i,j} \, \dm(v) \right)}_{i,j},\]

where ${F(v)}_{i,j}$ is the entry  $(i,j)$ of the matrix $F(v)$.

We will list a sequence of trivial results (without proof) that will be used next.

\begin{lemma}
	For an integrable transformation  $F: M_k \to M_k$
	\[\tr \int_{M_k} F(v) \, \dm(v) = \int_{M_k} \tr F(v) \, \dm(v).\]
\end{lemma}

\begin{lemma}
	Given a matrix $B\in M_k$ and an integrable transformation $F:M_k\to M_k$, then,
	\[ B\int_{M_k} F(v) \, \dm(v) = \int_{M_k} B F(v) \, \dm(v). \]
\end{lemma}

\begin{proposition}\label{fun_int}
	If $l: M_k \to \C$ is a linear functional and $F:M_k \to M_k$ is integrable, then,

	\[l\left(\int_{M_k} F(v) \, \dm(v) \right) = \int_{M_k} l(F(v)) \, \dm(v).\]
\end{proposition}

\begin{defi}\label{kui}
Given a measure $\mu$ on $M_k$ and  a measurable funtion $L:M_k\to M_k$, we say that $\mu$ is  $L$-square integrable, if

	\[\int_{M_k} \norm{L(v)}^2 \,\dm(v) < \infty.\]

For a fixed $L$ we denote by $\mathcal{M}(L)$ the set of  $L$-square-integrable measures. We also denote $\mathcal{P}(L)$ the set of  $L$-square-integrable probabilities.
\end{defi}

\begin{defi}\label{kiku}
  Given a measurable function  $L:M_k\to M_k$ and a  $L$-square-integrable measure $\mu$  we define the linear operator $\phi_L:M_k \to M_k$ via the expression
	\[\rho \,\to\,\phi_L(\rho) = \int_{M_k} L(v) \rho {L(v)}^\dagger \, \dm(v).\]
\end{defi}
\medskip

For  a given  $H: M_k \to M_k$ (which plays the role of a Hamiltonian)  we present a version of the Ruelle Theorem: a variational principle  of pressure related to an eigenvalue problem for a kind of Ruelle operator (see Theorem \ref{peqlog}).

Remember that if $A,B \in M_k$ with $A,B\ge 0$, then $\tr(AB) \le \tr(A) \tr(B)$. Therefore, if $\rho\ge 0$, we have
\begin{align*}
  \norm{\phi_L(\rho)}^2&=\tr(\phi_L(\rho){\phi_L(\rho)}^{\dagger}) \\
	&= \int_{M_k}\int_{M_k} \tr(L(v)\rho {L(v)}^{\dagger}{L(w)}^{\dagger}\rho L(w))\dm(v)\dm(w)\\
	&=\int_{M_k}\int_{M_k} \tr(\rho {L(v)}^{\dagger}{L(w)}^{\dagger}\rho  L(w)L(v))\,\,\dm(v)\dm(w) \\
	&\le \tr(\rho)  {\int}_{M_k}\int_{M_k} \tr(\rho L(w)L(v){L(v)}^{\dagger}{L(w)}^{\dagger})\,\,\dm(v)\dm(w) \\
	&\le {\tr(\rho)}^2 {\int}_{M_k}\int_{M_k} \tr(L(w)L(v){L(v)}^{\dagger}{L(w)}^{\dagger})\,\,\dm(v)\dm(w)\\
	&\le {\tr(\rho)}^2 {\int}_{M_k} \norm{L(v)}^2 \,\dm(v) \, \int_{M_k}\norm{L(w)}^2 \,\dm(w) < \infty.
\end{align*}

For a general $\rho \in M_k$, we write $\rho = \rho_{+} - \rho_{-}$ where $\rho_{+} = \abs{\rho}$ and $\rho_{-} = \abs{\rho} - \rho$ are both positive semidefinite matrices. By linearity of $\phi_L$, we have

\[\phi_L(\rho) = \phi_L(\rho_{+}) - \phi_L(\rho_{-}),\]

\noindent hence, $\phi_L$ is well defined.

\medskip

\begin{proposition}\label{fun_int1}
	Given a measurable function $L:M_k \to M_k$ and a $L$-square integrable  measure $\mu$, then,  the dual transformation $\phi_L^{*}$  is given by
	\[\phi^*_L(\rho) = \int_{M_k} {L(v)}^{\dagger} \rho L(v) \, \dm(v).\]
\end{proposition}

\begin{defi}
Given a measurable function $L:M_k \to M_k$ and a $L$-square integrable
measure $\mu$ over $M_{k}$, then, the transformation $\phi_{L}$ is called
\textbf{stochastic} if
	\[\phi_L^*(\Id_k) = \int_{M_k} {L(v)}^{\dagger} L(v) \, \dm(v) = \Id_k.\]

	By abuse of language, we sometimes say $L$ stochastic to mean that $\phi_{L}$ is stochastic.
\end{defi}

We will be able to define the concept of entropy when the $\phi_L$ is stochastic.

\begin{defi}\label{posu}
  A linear map $\phi: M_k \to M_k$ is called {\bf positive} if takes positive matrices to positive matrices.
\end{defi}

\begin{defi}\label{posu1}
  A positive linear map $\phi: M_k \to M_k$ is called {\bf completely positive}, if for any $m$, the linear map $\phi_m=\phi \otimes I_m: M_k \otimes M_m \to M_k \otimes M_m$ is positive, where $I_m$ is the identity operator acting on the matrices in $M_m$.
\end{defi}

\begin{defi} If $\phi:M_k\to M_k$ is a linear map and satisfies
\begin{enumerate}
		\item $\phi$ is completely positive;
		\item $\phi$ preserves trace.
	\end{enumerate}

\noindent	Then, we say that $\phi$ is a quantum channel.
\end{defi}

\begin{theorem}\label{quantum_channel_phi}

  Given $L:M_{k}\to M_{k}$ and $\mu$ a $L$-square measure. Then the associated transformation $\phi_L$  is completely positive. Moreover, if  $\phi_L$ is stochastic then preserves trace.
\end{theorem}
\begin{proof}
\textit{1. $\phi_L$ is completely positive:} suppose $A \otimes B\in M_n\otimes M_k$ satisfies $A\otimes B \ge 0$ and $\psi \in \C^n\otimes \C^k$. Then, if $\psi_L(v) = (Id_n\otimes {L(v)}^{\dagger})\psi$ we get

\begin{align*}
	\bra{\psi} A\otimes \phi_L(B) \ket{\psi} &= \bra{\psi} A\otimes \int_{M_k} L(v) B {L(v)}^\dagger \, \dm(v) \ket{\psi}	\\
		&= \int_{M_k} \bra{\psi} A \otimes ( L(v) B {L(v)}^\dagger ) \ket{\psi} \, \dm(v) \\
		&= \int_{M_k} \bra{\psi} (\Id_n\otimes L(v))(A \otimes B)(Id_n\otimes {L(v)}^\dagger) \ket{\psi} \, \dm(v)\\
		&= \int_{M_k} \bra{(\Id_{n}\otimes {L(v)}^\dagger)\psi}(A \otimes B) \ket{(Id_n\otimes {L(v)}^\dagger)\psi} \, \dm(v)\\
		&= \int_{M_k} \bra{\psi_L(v)}(A \otimes B) \ket{\psi_L(v)} \, \dm(v) \ge 0.
\end{align*}

Above we use the positivity of  $A\otimes B$ in order to get $\bra{\psi_L(v)}(A \otimes B) \ket{\psi_L(v)} \ge 0$. We also used in some of the equalities  the fact that  $l(X) := \bra{\psi} A\otimes X \ket{\psi}$  is a linear functional and therefore we can apply proposition $\ref{fun_int}$.

\textit{2. Under our assumption $\phi_L$ preserves trace:} given $B\in M_k$
\begin{align*}
\tr \phi_L(B) &= \tr \left(\int_{M_k} L(v) B {L(v)}^\dagger \, \dm(v)\right)\\
&= \int_{M_k} \tr \left( L(v) B {L(v)}^\dagger \right) \, \dm(v)\\
&= \int_{M_k} \tr \left( B {L(v)}^\dagger L(v) \right) \,\dm(v)\\
&= \tr \left(B \int_{M_k} {L(v)}^\dagger L(v) \,\dm(v) \right) \\
&= \tr(B).
\end{align*}
\end{proof}

\begin{remark}[$\phi_L^*$ is completely positive]
When $L$ is measurable, then, using the same reasoning as above one can show that $\phi_L^*$  is completely positive.

We say that $\phi_L$ preserves unity if $\phi_L(\Id) = \Id$. In this case, $\phi_L^*$ preserves trace. When $\phi_L^*$ preserves the identity then $\phi_L$ preserves trace.
\end{remark}

\section{The eigenvalue property for $\phi_L$}\label{eig}

In  this section, we will investigate questions related to the existence of eigenvalues and eigenmatrices for the setting of  Quantum Information. An important role will be played by a result about  positive maps on $C^*$-algebras described in  \cite{evans1978spectral} which presents a noncommutative version of the Perron Theorem.

\begin{defi}[Irreducibility]\label{gilo}
	We say that $\phi: M_k \to M_k$ is irreducible if one of the equivalent properties is true
	\begin{itemize}
		\item Does not exists $\lambda > 0$ and a projection $p$ on a proper non-trivial subspace of $\C^k$, such that, $\phi(p) \le \lambda p$;
		\item For all non null $A\ge 0$, ${(\Id+\phi)}^{k-1}(A) > 0$;
		\item For all non null $A\ge 0$ there exists $t_A>0$, such that, $(e^{t_A \phi})(A) > 0$;
		\item If $P\in M_{k}$ is a hermitian projector such that $\phi(PM_{k}P)\subset PM_{k}P$, then $P\in\{0,\Id\}$;
		\item For all pair of non null positive matrices  $A,B\in M_k$ there exists a natural number $n\in \{1,\ldots,k-1\}$, such that, $\tr[B\phi^n(A)] > 0$.
	\end{itemize}
\end{defi}

The proof of the equivalence of the two first items appears in $\cite{evans1978spectral}$.

The equivalence of the two middle ones appears in $\cite{schrader2000perron}$ where also one can find the proof of the improved positivity (to be defined below) which implies irreducibility. For the proof that the last two items, we refer to $\cite{wolf2012quantum}$.

\begin{defi}[Irreducibility]\label{gil}
	Given $\mu$ we will say (by abuse of language) that $L$ is irreducible for $\mu$ or $\mu$-irreducible if the associated $\phi_L$ is irreducible.
\end{defi}

\begin{lemma}\label{irred_E_Ck_equiv}
  Given $L: M_{k}\to M_{k}$ and $\mu$ a $L$-square measure, the following statements are equivalent:
  \begin{enumerate}
    \item $\phi_{L}$ is irreducible;
	\item If $E\subset\C^{k}$ is a
	  subspace such that $L(v)E \subset E$ for all $v\in\supp\mu$, then $E\in\{\, \{0\},\C^{k}\}$.
  \end{enumerate}
\end{lemma}

\begin{proof}
  $1. \to 2$.: If $\phi_{L}$ is irreducible and $E\subset\C^{k}$ is a
  subspace such that $L(v)E \subset E$ for all $v\in\supp\mu$, take $P$ the
  orthogonal projection on $E$. Then $PL(v)P = L(v)P$ for all $v\in\supp\mu$.
  Moreover, for every $A\in M_{k}$

  \begin{align*}
	\phi_{L}(PAP) &= \int_{M_{k}} L(v)PAP{L(v)}^{\dagger}\dm(v)\\
	   	&= \int_{\supp\mu} PL(v)PAP{L(v)}^{\dagger}P\dm(v)\\
	   	&= P\int_{\supp\mu}L(v)PAP{L(v)}^{\dagger}\dm(v)P\in PM_{k}P,
  \end{align*}

\noindent and by the fourth equivalence of $\ref{gilo}$, $P\in\{0,\Id\}$. Therefore $E=\{0\}$
  or $E=\mathbb{C}^{k}$.

  $2.\to1$.: If there is $P \in M_{k}$ Hermitian projection such that
  $\phi_{L}(PM_{k}P)\in PM_{k}P$, take $E=\text{\normalfont{Im}}\,P$,
  $x\in E$ and $A=\ket{x}\bra{x}$. Then we have

  \begin{align*}
      0 &= \tr(\phi_{L}(PAP) - P\phi_{L}(PAP)P) =\\
        &= \int_{M_{k}} \tr\left(L(v)A{L(v)}^{\dagger} - PL(v)A{L(v)}^{\dagger}P\right) \dm(v)\\
        &= \int_{B_{x}} \tr\left(L(v)A{L(v)}^{\dagger} - PL(v)A{L(v)}^{\dagger}P\right) \dm(v)\\
        &= \int_{B_{x}} \tr\left(L(v)A{L(v)}^{\dagger}(I-P)\right) \dm(v),
  \end{align*}

\noindent  where $B_{x}:=\{v\,\vert\,L(v)x\notin E\}$. Suppose $P\notin\{0, \Id\}$,
    then we have
  \[\tr(L(v)A{L(v)}^{\dagger}) - \tr(PL(v)A{L(v)}^{\dagger}P) > 0,\]
  and, since the integral is zero, $\mu(B_{x}) = 0$.
  Thus, for all $x\in E$, $\supp\mu \subset B_{x}^{c}$ and so $L(v)E \subset E$
  for $v \in \supp\mu$. By hypothesis, we have that $E\in\{\, \{0\}, \C^{k}\}$
  which brings us to an absurd.

\end{proof}

\begin{defi}[Improving positivity]
	We say that  $\phi_L$ is positivity improving, if $\phi_L(A)>0$, for any non-null $A\ge 0$.	Note that improving positivity implies irreducibility.
\end{defi}

	For any  $\mu$ and square-integrable $L$ the Theorem $\ref{quantum_channel_phi}$ assures that $\phi_L$ is completely positive. In the case $\phi_L$ is irreducible we can use the Theorem 2.3 and 2.4 of $\cite{evans1978spectral}$ in order to get  $\lambda$ and $\rho > 0$, such that, $\phi_L(\rho) = \lambda \rho$ ($\rho$ is unique up to multiplication by scalar). For what comes next, we will choose $\rho$ such that $\tr \rho = 1$. Moreover, in the same work the authors show that $\phi_L$  is irreducible, if and only if, $\phi_L^*$ also is completely positive, and therefore we get:

\begin{theorem}[The spectral radius is a simple eigenvalue]\label{raio_espec}
	Given a square integrable $L:M_k\to M_k$ assume that  the associated $\phi_L$ is irreducible. On a Hilbert space, the spectral radius $\lambda_L> 0$ of $\phi_L$ and $\phi_L^*$ is the same. In this case it is also an eigenvalue and it is simple. We denote, respectively, by $\rho_L > 0$ and $\sigma_L > 0$, the eigenmatrices, such that,  $\phi_L(\rho_L) = \lambda_L\rho_L$ and $\phi_L^*(\sigma_L) = \lambda_L\sigma_L$, where $\rho_L$ and $\sigma_L$ are the unique non null eigenmatrices  (up to multiplication by scalar).
\end{theorem}

	The above theorem is the natural  version of the Perron-Frobenius Theorem for the present setting. It is natural to think that $\phi_L$ acts on density matrices and $\phi_L^*$ acts in selfadjoint matrices.
\begin{remark}\label{unique}
    We choose $\rho_L$ in such way that  $\tr \rho_L = 1$ and after that, we take $\sigma_L$ such that $\tr (\sigma_L \rho_L\, )= 1$. By doing that, we have chosen the precise scalar multiples that makes both $\rho_L$ and $\sigma_L \rho_L$ densities. Notice that, as eigendensity, $\rho_L$ is unique.
	We point out that at this moment it is natural to make an analogy with Thermodynamic Formalism: $\phi_L^*$ corresponds to the Ruelle operator (acting on functions) and $\phi_L$ to the dual of the Ruelle operator (acting on probabilities). We refer the reader to $\cite{PP}$ for details. In this sense, the density operator $\sigma_L \rho_L$ plays the role of an equilibrium probability. The paper $\cite{Sp}$ by Spitzer describes this formalism in a simple way in the case the potential depends on two coordinates.
\end{remark}

\begin{remark}
  If $L$ is irreducible and stochastic (resp. $\phi_{L}$ is unital, i.e., $\phi_{L}(\Id) = \Id$) then
  $\lambda_{L}=1$ and $\sigma_{L}=\Id_{k}$ (resp. $\rho_{L}=\Id_{k}$) by Proposition 6.1 on $\cite{wolf2012quantum}$ page 91.
\end{remark}

\subsection{Normalization}

We consider in this section a fixed measure $\mu$  over $M_k$ which plays the role of the \textit{a priori} probability.

In this section we will introduce the concept of normalized transformation $L$ (see definition \ref{rert}). If $L$ is not normalized we will be able to find an
associated $\hat{L}$ which is normalized (see \eqref{opii}).

Given a continuous $L$ {\bf (variable)} we assume in this section that  $\phi_L$ is irreducible  (we do not assume that preserves trace).

We will associate to this square integrable transformation $L:M_k \to M_k$ (and the associated $\phi_L$) another transformation $\hat{L}: M_k \to M_k$  which will correspond to a normalization of $L$. This will define another quantum channel $\phi_{\hat{L} }: M_k \to M_k$.

Results of this section have a large intersection with some material  in $\cite{wolf2012quantum}$. For completeness, we describe here what we will need later.

Consider $\sigma_L$ e $\lambda_L$ as described above.
 As $\sigma_L$ is positive we consider ${\sigma_L}^{1/2} > 0$ and ${\sigma_L}^{-1/2} > 0$.

 In this way we define 
 
 \begin{equation} \label{opii}\hat{L}(v) =\frac{1}{\sqrt{\lambda_L}} {\sigma_L}^{1/2} L(v) {\sigma_L}^{-1/2}.
 \end{equation}

Using the measure $\mu$ we can define  the associated $\phi_{\hat{L}}$.

 Therefore,
			\begin{align*}
				\phi_{\hat{L}}^*(\Id) &=\frac{1}{\lambda_L} \int_{M_k} {\sigma_L}^{-1/2} {L(v)}^\dagger {\sigma_L}^{1/2} {\sigma_L}^{1/2} L(v) \sigma_L^{-1/2}\dm\\
				&= \frac{1}{\lambda_L}\sigma_L^{-1/2}  \int_{M_k} {L(v)}^\dagger \sigma_L L(v) \dm \,\,\sigma_L^{-1/2}\\
				&= \frac{1}{\lambda_L}\sigma_L^{-1/2} \phi_L^*(\sigma_L) \,\sigma_L^{-1/2} \\
				&= \frac{1}{\lambda_L}\sigma_L^{-1/2} \lambda_L\sigma_L \,\sigma_L^{-1/2}\\
				&= \sigma_L^{-1/2} \sigma_L  \sigma_L^{-1/2}\\
				&= \Id.
			\end{align*}

			Note that ${\hat{L}(v)}^\dagger = \frac{1}{\sqrt{\lambda_L}} \sigma_L^{-1/2} {L(v)}^\dagger \sigma_L^{1/2}$. From this we get easily  that  $\phi_{\hat{L}}$ is completely positive and preserves trace (is stochastic).

We will show  that $\phi_{\hat{L}}$ is irreducible. Given $A\in M_k$ we have

				\[\phi_{\hat{L}}(A) = \frac{1}{\lambda_L}\sigma_L^{1/2} \phi_L(\sigma_L^{-1/2}A\sigma_L^{-1/2}) \sigma_L^{1/2}.\]

			Then,
			\begin{align*}
				\phi_{\hat{L}}^2(A) &= \frac{1}{\lambda_L}\sigma_L^{1/2} \phi_L(\sigma_L^{-1/2}\frac{1}{\lambda_L}\sigma_L^{1/2} \phi_L(\sigma_L^{-1/2}A\sigma_L^{-1/2}) \sigma_L^{1/2}\sigma_L^{-1/2}) \sigma_L^{1/2}\\
				&=\frac{1}{\lambda_L^2}\sigma_L^{1/2} \phi_L^2(\sigma_L^{-1/2}A\sigma_L^{-1/2}) \sigma_L^{1/2}.
			\end{align*}

			By induction we get
				\[\phi_{\hat{L}}^n(A) = \frac{1}{\lambda_L^n}\sigma_L^{1/2} \phi_L^n(\sigma_L^{-1/2}A\sigma_L^{-1/2}) \sigma_L^{1/2}.\]

			Given $A,B\ge 0$, note that $\sigma_L^{-1/2}A\sigma_L^{-1/2}\ge 0$ and $\sigma_L^{1/2}B\sigma_L^{1/2}\ge 0$. Therefore, using irreducibility of $\phi_L$, there exists an integer $n\in \{1,\ldots, k-1\}$, such that,
			\begin{align*}
				0 &< \lambda_L^{-n}\tr[\sigma_L^{1/2}B\sigma_L^{1/2}\phi_L^n(\sigma_L^{-1/2}A\sigma_L^{-1/2})]\\
				&=\lambda_L^{-n}\tr[B\sigma_L^{1/2}\phi_L^n(\sigma_L^{-1/2}A\sigma_L^{-1/2})\sigma_L^{1/2}]\\
				&=\tr[B\phi_{\hat{L}}^n(A)].
			\end{align*}

			Therefore, $\phi_{\hat{L}}$ is irreducible and completely positive and preserves trace. In this way, to the given  $L$ we can associate $\hat{L}$ which will be called the normalization of $L$. The transformation $\phi_{\hat{L}}$ is a quantum channel.

\begin{defi}  Given the measure $\mu$  over $M_k$ we denote by $\mathfrak{L}(\mu)$ the set of all integrable $L$ such that  the associated  $\phi_L$ is irreducible.

\end{defi}

\begin{defi} \label{rert}
	Suppose   $L$ is in $\mathfrak{L}(\mu)$. We say that $L$ is normalized if $\phi_L$ has spectral radius $1$ and preserves trace. We denote by $\mathfrak{N}(\mu)$ the set of all normalized $L$.
\end{defi}

Note that the transformation $\hat{L}$ defined above in \eqref{opii} is normalized.

If $L \in  \mathfrak{N}(\mu)$, then, we get from  Theorem $\ref{raio_espec}$ and the fact that  $\phi_L^*(\Id) = \Id$, that the spectral radius, which is also a simple eigenvalue, is $\lambda_L = 1$. According to Remark $\ref{unique}$, there is a unique eigendensity $\rho_L$ such that $\phi_L(\rho_L) = \rho_L$. These properties will be important for what will come next.

\begin{theorem}[Ergodicity and temporal means]
Suppose $L \in \mathfrak{N}(\mu)$. Then, for all density matrix $\rho \in M_k$ it is true that
\[ \lim_{N \to \infty} \frac{1}{N}\sum_{n=1}^{N} \phi_L^n(\rho) = \rho_L,\]
where $\rho_L$ is the density matrix associated to $L$.
\end{theorem}
\begin{proof}
	The proof follows from  Theorem $\ref{raio_espec}$ and Corollary  6.3 in $\cite{wolf2012quantum}$.
\end{proof}


The above result connects irreducibility and ergodicity (the temporal means have a unique limit).

\section{Entropy}\label{ent}

	In this section, we will define entropy for $\phi_L$, when the associated $L$   is irreducible and stochastic (see Definition \ref{oou1}). After that, it will be possible to give a meaning for a certain   variational principle of pressure in Definition \ref{oou} (this is similar to the setting in Thermodynamic Formalism which is described in $\cite{PP}$, for instance).

Remember the classical entropy is defined just for invariant (stationary) probabilities. Something of this sort is required for defining the entropy of a quantum channel $\phi_L$: $L$ has to be stochastic. These $\phi_L$ will play in some sense the role of the different possible invariant probabilities.

We will explore some ideas which were already present on the paper $\cite{baraviera2010thermodynamic}$ (which explores some previous nice results on $\cite{lozinski2003quantum}$ and $\cite{slomczynski2003dynamical}$) which considers a certain \textit{a priori} probability.

Hereafter, we consider {\bf fixed a measure $\mu$}  over $M_k$ which plays the role of the \textit{a priori} probability.  Given $L\in \mathfrak{L}(\mu)$ we will associate in a natural way the transformation $\phi_L: M_k \to M_k$.

\begin{defi}
We denote by $\phi=\phi_\mu$ the set of all $L$ such that the associated
$\phi_L: M_k \to M_k$ is irreducible and stochastic.
\end{defi}

				We will describe a discrete-time  process that takes values on $M_k$.

		Suppose  $L$ is irreducible and stochastic.
		We  will associate to such $L$ a kind of ``transition probability kernel'' $ P_L$ (to be defined soon) acting on matrices. Given the matrices $v$ and $w$ the value $P_L(v,w)$ will describe the probability of going in the next step  to $w$ if the process is on $v$.

Given $L$, suppose that the discrete-time process is given in such a way that the initial state is  described by the density matrix $\rho_L$ which is invariant for $\phi_L$ (see Theorem $\ref{raio_espec}$).

The reasoning here is that such process should be in ``some sense stationary'' because $\rho_L$
is invariant by $\phi_L$. As we said before in ergodic theory the concept of Shanon-Kolmogorov entropy has a meaning just for invariant   (for a discrete-time dynamical system) probabilities. Therefore, something of this order is required.

 In our reasoning given that the state is described by  $\rho$, then,  in the next step of the process we get $\frac{L(v)\rho {L(v)}^{\dagger}}{\tr(L(v)\rho {L(v)}^{\dagger})}$ with probability $\tr(L(v)\rho {L(v)}^{\dagger})\,\dm(v)$.

 This discrete-time process takes values on density operators in  $M_k$.

 \begin{defi} \label{oou1}
 We define entropy for $L$ (or, for $\phi_L$) by the expression (when finite):

	\[h(L)= h_\mu(L) := - \int_{M_k \times M_k\hspace{-15px}} \tr(L(v)\rho_L {L(v)}^{\dagger}) P_L(v,w)\log P_L(v,w) \,\,\dm(v)\dm(w),\]

		where
		\[ P_L(v,w) := \frac{\tr(L(w)L(v)\rho_L {L(v)}^{\dagger}{L(w)}^{\dagger})}{\tr(L(v)\rho_L {L(v)}^{\dagger})}.\]
		\end{defi}

		This definition is a generalization of the analogous concept presented on the papers $\cite{baraviera2010thermodynamic}$, $\cite{BLLT2}$ and $\cite{BLLT1}$.

		Note that $\tr(L(v)\rho_L {L(v)}^{\dagger})$ is the probability of being in  state  $\frac{L(v)\rho_L {L(v)}^{\dagger}}{\tr(L(v)\rho_L {L(v)}^{\dagger})}$.
		Moreover, $P_L(v,w)$ describes the probability of going from $v$ to $w$, being in state\\
		\[\frac{L(w)L(v)\rho_L {L(v)}^{\dagger}{L(w)}^{\dagger}}{\tr(L(w)L(v)\rho_L {L(v)}^{\dagger}{L(w)}^{\dagger})}.\]
		In this way  $h_\mu(L)$ in some way     resembles the analogous expression of entropy for the case of Markov chains.

		We will show in Example $\ref{Mark1}$ that the above definition of entropy is indeed a natural generalization of the classical one in Ergodic Theory.

		Suppose $H: M_k \to M_k$ is square integrable, irreducible and $H(v)\ne 0$, for $\mu$\emph{-a.e.} $v$. For such $H$,  consider  the  corresponding $\rho_H$, $\sigma_H$ and $\lambda_H$ which are given by Theorem $\ref{raio_espec}$,  where $\tr \rho_H = 1$ and $\tr \sigma_H \rho_H = 1$.

		This $H$ describes the action of a potential.

		 Then, we define

		\[U_H(v) := \log\left(\tr(\sigma_H H(v)\rho_H {H(v)}^{\dagger}) \right).\]

		\begin{defi} \label{oou}
		We define the pressure of $H$ by
		\[P_\mu(H)=P(H) := \sup_{L \in \phi}\left\{ h_\mu(L) + \int U_H(v) \,\, \tr(L(v)\rho_L {L(v)}^{\dagger}) \,\, \dm(v) \right\}.\]
		\end{defi}

		Remember that $\phi_\mu$ is the set of all  $L: M_k \to M_k$ which are square-integrable, irreducible, and stochastic.

		\begin{defi}\label{ort}
		  Given $\mu$ and $H$ as above we say that $\phi_L$, for some $L \in \phi_\mu$, is a Gibbs channel, if
		\[P_\mu(H)=  h_{\mu}(L) + \int U_H(v) \,\, \tr(L(v)\rho_L {L(v)}^{\dagger})\dm(v).\]
		\end{defi}

We will need soon the following well-known result (see $\cite{PP}$).

		\begin{proposition}\label{gibbs_inequality}
			Suppose  $p,q:M_k \to \R_{+}$ are such that  $p,q > 0$, $\mu$-almost everywhere, $\int_{M_k} p \,\,\dm = 1$ and $\int_{M_k} q \,\,\dm = 1$. Then,
				\[- \int p \log p \,\, \dm + \int p \log q \,\, \dm \le 0.\]

			Moreover, the above inequality is an  equality just when $p = q$, $\mu$-almost everywhere.
		\end{proposition}

		\begin{theorem}\label{var_prin_1}
			Assume that $H: M_k \to M_k$ is continuous,  irreducible and $H(v)\neq 0$ for $\mu$-a.e. $v$, then,
				\[P(H) := \sup_{L\in \phi}\left\{ h_\mu(L) + \int U_H(v) \,\, \tr(L(v)\rho_L {L(v)}^{\dagger}) \,\, \dm(v) \right\} \le \log(\lambda_H), \]

			The  supremum is attained only if
				\[\frac{\tr(L(w)L(v)\rho_L {L(v)}^{\dagger}{L(w)}^{\dagger})}{\tr(L(v)\rho_L {L(v)}^{\dagger})} = \frac{1}{\lambda_H}\tr(\sigma_H H(w)\rho_H H(w)^{\dagger}),\quad\text{for}\,\mu\text{-a.e.}\, v,w.\]

			In this case, $P(H) = \log(\lambda_{H})$.
		\end{theorem}

\medskip

\begin{proof}
		We define $q(w) :=  \frac{1}{\lambda_H}\tr(\sigma_H H(w)\rho_H H(w)^{\dagger})$. Note that
		\begin{align*}
			\int q\,\,\dm &= \frac{1}{\lambda_H}\int \tr(\sigma_H H(w)\rho_H H(w)^{\dagger})\,\,\dm(w)\\
			&= \frac{1}{\lambda_H}\tr\left(\sigma_H\int H(w)\rho_H H(v)^{\dagger}\, \dm(w)\right)\,\,\\
			&= \frac{1}{\lambda_H} \tr(\sigma_H \lambda_H\rho_H)\\
			&= \tr(\sigma_H\rho_H) = 1.
		\end{align*}

		For fixed $v$ and irreducible and stochastic $L$ take
			$$p_v(w) = P_L(v,w) = \frac{\tr(L(w)L(v)\rho_L L(v)^{\dagger}L(w)^{\dagger})}{\tr(L(v)\rho_L L(v)^{\dagger})},$$

		If $\tr(L(v)\rho_L L(v)^{\dagger})\ne0$ and $p_{v}(w) = 0$ otherwise.
		It follows that

		\begin{align*}
			\int p_v(w) \,\,\dm(w) &= \int \frac{\tr(L(w)L(v)\rho_L L(v)^{\dagger}L(w)^{\dagger})}{\tr(L(v)\rho_L L(v)^{\dagger})}\,\, \dm(w)\\
			 &=\frac{\tr(L(v)\rho_L L(v)^{\dagger}\int L(w)^{\dagger}L(w))\,\,\dm(w)}{\tr(L(v)\rho_L L(v)^{\dagger})} = 1.
		\end{align*}

		From Proposition $\ref{gibbs_inequality}$ we get that for each  $v$
		\begin{equation}\label{pq_gibbs_inequality}
			- \int p_v(w)\log(p_v(w))\,\,\dm(w) + \int p_v(w) \log (q(w))\,\,\dm(w) \le 0.
		\end{equation}

		Equality will happen when

		$$ \frac{\tr(L(w)L(v)\rho_L L(v)^{\dagger}L(w)^{\dagger})}{\tr(L(v)\rho_L L(v)^{\dagger})} = \frac{1}{\lambda_H}\tr(\sigma_H H(w)\rho_H H(w)^{\dagger}),$$

		for $\mu$-almost everywhere   $w$.

		Note that from ($\ref{pq_gibbs_inequality}$) it follows that

		\begin{align*}
			\int -P_L(v,w) \log P_L(v,w) &+ P_L(v,w) \log\left(\tr(\sigma_H H(w)\rho_H H(w)^{\dagger})\right) \dm(w) \\
			&\le \int P_L(v,w) \log(\lambda_H)\,\,\dm(w) = \log(\lambda_H).
		\end{align*}

Now we multiply  both sides of the above inequality by $\tr(L(v)\rho_L {L(v)}^{\dagger})$,
integrate with respect to $v$ (remember that $\int\tr(L(v))\rho_L {L(v)}^{\dagger}) = 1$)
and we get

\begin{align*}
	h_\mu&(L)\,+\\ &\int \tr(L(w)L(v)\rho_L L(v)^{\dagger}L(w)^{\dagger})\log\left(\tr(\sigma_H H(w)\rho_H H(w)^{\dagger})\right)\,\,\dm(w)\dm(v)\\
	&= h_\mu(L) + \int \tr(L(w)\phi_L(\rho_L) L(w)^{\dagger})\log\left(\tr(\sigma_H H(w)\rho_H H(w)^{\dagger})\right)\,\,\dm(w)\\
	&= h_\mu(L) + \int \log\left(\tr(\sigma_H H(v)\rho_H H(v)^{\dagger})\right)\,\tr(L(v)\rho_L L(v)^{\dagger})\,\,\dm(v)\\
	&\le \log(\lambda_H).
\end{align*}

As this is true for any $L\in \phi$, we take  the $\sup$ over all such $L$ to finally get:

	$$P(H) \le \log(\lambda_H).$$

\end{proof}
\medskip

A natural question: is there a $L\in \phi$ such that the supremum is attained? This kind of result would correspond in our setting to the Ruelle Theorem of Thermodynamic Formalism (see \cite{PP}). In this direction, we are able to get Theorem \ref{peqlog}.

\medskip

Before trying to address this question  we point out that given $H$ as above one can get the associated normalized  $\hat{H}$ by the expression       $\hat{H}= \frac{1}{\sqrt{\lambda_H}}\sigma_H^{1/2} H \sigma_H^{-1/2}$.

	Note that $\sigma_{\hat{H}} = \Id$, $\rho_{\hat{H}} = \sigma_H^{1/2}\rho_H\sigma_H^{1/2}$ and $\lambda_{\hat{H}} = 1$. Therefore,

\begin{align*}
	\int& \log\left(\tr(\sigma_{\hat{H}} \hat{H}(v)\rho_{\hat{H}} \hat{H}(v)^{\dagger})\right)\,\tr(L(v)\rho_L L(v)^{\dagger})\,\,\dm(v)\\
	 &=\int \log\left(\tr(\frac{1}{\lambda_H}\sigma_H^{1/2} H(v) \sigma_H^{-1/2}\sigma_H^{1/2}\rho_H\sigma_H^{1/2}\sigma_H^{-1/2} H(v)^{\dagger} \sigma_H^{1/2} )\right)\,
      \tr(L(v)\rho_L L(v)^{\dagger})\,\,\dm(v)\\
	&=\int \log\left(\tr( \frac{1}{\lambda_H}\sigma_H H(v)\rho_H H(v)^{\dagger} )\right)\,\tr(L(v)\rho_L L(v)^{\dagger})\,\,\dm(v)\\
	&=\int \log\left(\tr(\sigma_H H(v)\rho_H H(v)^{\dagger} )\right)\,\tr(L(v)\rho_L L(v)^{\dagger})\,\,\dm(v) - \log(\lambda_H).
\end{align*}

From the above reasoning we get:

\begin{theorem}\label{var_prin_2}
Assume that $H: M_k \to M_k$ is irreducible, square integrable and $H(v)\neq 0$, for $\mu$-a.e. $v$. If $\hat{H}$ denotes the associated normalization, then,
	$$P(\hat{H}) = P(H) - \log(\lambda_H).$$
\end{theorem}

Note that $\hat{H}\in \phi_\mu$.

\begin{theorem}\label{peqlog}
	If $H$ is irreducible, square integrable and $H(v)\neq 0$, for $\mu$-a.e. $v$, then,
	$$P(H) = \log \lambda_H.$$
\end{theorem}
\begin{proof}
We already know that $P(H)\le \log \lambda_H$. We will show that  there exists an irreducible and stochastic $L$ which attains the supremum. In order to do that  we take an orthonormal  basis $\{\ket{i}\}_{i=1,2,\ldots,k}$ of $\C^k$. Then, we define an  operator  $P$ such that $P\ket{i+1} = \ket{i}$ (for instance, $P = \sum_{i=1}^k \ket{i}\bra{i+1}$ and by convention $\ket{1}=\ket{k+1}$).

		Note that the dual of $P$ is $P^\dagger = \sum_i \ket{i+1}\bra{i}$. This is so because given $u,v\in \C^k$, we get that

			$$\inner{u}{Pv} = \sum_i \inner{u}{\ket{i}\bra{i+1} v} = \sum_i \inner{\ket{i+1}\bra{i}u}{v} = \inner{P^{\dagger}u}{v}.$$

		Moreover, $P^{\dagger}P = \Id$. Indeed,

			$$\sum_{i,j}  \ket{j+1}\bra{j}\ket{i}\bra{i+1} = \sum_i \ket{i}\bra{i} = \Id.$$

		Now, take $Q=(q_{ij})$ the matrix with $q_{kk} = -1$, $q_{ii}=1$, for $i=1,...,k-1$, and $q_{ij}=0$ otherwise. Note that $Q^\dagger Q = \Id$.\\

		Consider $\rho_H, \sigma_H, \lambda_H$ given by Theorem $\ref{raio_espec}$, where $\tr(\sigma_H\rho_H)=1$ and $\tr(\rho_H) = 1$ and
		let $\varphi(v) = \sqrt{\frac{1}{\lambda_H}\tr(\sigma_H H(v) \rho_H {H(v)}^{\dagger})}$.

		Note that if $\#\supp\mu = 1$, $H$ can't be irreducible because any eingenvector of $H(v)$ for $v\in\supp\mu$ generates an invariant subspace.

		There exist $v_1,v_2 \in \supp\mu$ with $\varphi(v_i) \ne 0$ by hypothesis. Take $O$ an open set with $v_1\in O$ and $d(\overline{O}, v_2) > 0$. Now we can define $L$ by $L(v) = \varphi(v)P$, for $v \notin O$, and $L(v) = \varphi(v)Q$, for $v \in O$. \\

		Observe that $L(v)^\dagger L(v) = \abs{\varphi(v)}^2\Id$, for all $v$, and $\int \abs{\varphi(v)}^2 \,\dm(v) = 1$. This implies that $\phi^*_L(\Id) = \Id$. \\

        Suppose that $E$ is an invariant subspace of $\C^k$ for all $L(v)$ with $v \in \supp\mu$. Of course, as $\varphi(v_i)\ne 0$, $E$ is invariant for $P$ and $Q$. In this sense, taking $x=(x_1,\ldots,x_k)\in E$, we get $Qx = (x_1,\ldots,-x_k) \in E$. As $E$ is a linear subspace this implies that $(x_1,\ldots,x_{k-1},0) \in E$, and $(0, \ldots, 0, x_k) \in E$. Taking $P^n(0,\ldots,0,x_k)$, for $n=0,\ldots,k-1$, if $x_k \ne 0$, we get a base of $\C^k$ in $E$. Therefore, if $x_k \ne 0$, we have $E = \C^k$. On the other hand, if initially $x_k = 0$, we take $P^n x$, where $(P^nx)_k \ne 0$, and we  use the previous argument. If there is no $x \in E$ and $n$ such that $(P^nx)_k \ne 0$, then $E=\{0\}$. Therefore, $\phi_L$ is irreducible by Lemma $\ref{irred_E_Ck_equiv}$.\\

		To show that $L$ satisfy the supremum for pressure, from the inequality give by Theorem $\ref{gibbs_inequality}$, it is enough to show that

		$$\frac{\tr(L(w)L(v)\rho_L L(v)^{\dagger}L(w)^{\dagger})}{\tr(L(v)\rho_L L(v)^{\dagger})} = \frac{1}{\lambda_H}\tr(\sigma_H H(w)\rho_H H(w)^{\dagger}).$$

		\bigskip

		In order to get this, observe that
		\begin{align*}
		 	\tr&(L(w)L(v)\rho_L L(v)^{\dagger}L(w)^{\dagger})\\
			&=\tr(L(v)\rho_L L(v)^{\dagger}L(w)^{\dagger}L(w))\\
			&=\tr(L(v)\rho_L L(v)^{\dagger})\abs{\varphi(w)}^2.
		\end{align*}

		Thus, the required equation holds.
\end{proof}

\section{ Process $X_n$, $n\in \mathbb{N}$, taking values on $ P(\C^k)$ }\label{pro1}

Consider a fixed measure $\mu$ on $M_k$ and a fixed $L:M_k \to M_k$, such that, $\int_{M_k} \norm{L(v)}^2\,\dm(v)  < \infty$, and, also  that $\phi_L$ is irreducible and stochastic.

Note that if, for example, $\mu$ is a probability and the the function $v \to \norm{L(v)}$ is bounded we get that
$\int_{M_k} \norm{L(v)}^2\,\dm(v)  < \infty$.

Denote by $P(\C^k)$ the projective space on  $\C^k$ with the metric $d(\hat{x},\hat{y})	= (1-\abs{\inner{x}{y}}^2)^{1/2}$, where $x,y$   are representatives with norm $1$ and $\inner{\cdot}{\cdot}$ is the canonical inner product.

We choose representatives and from now on for generic $\hat{x}, \hat{y}$   the associated ones are denoted by $x,y$.
We assume ``continuity'' on these choices.

Take $\hat{x} \in P(\C^k)$ and $S\subset P(\C^k)$. For a stochastic  $\phi_L$ we consider the kernel

\begin{equation}\label{jajo}
  \Pi_L(\hat{x}, S) = \int_{M_k} \textbf{1}_{S}(L(v)\cdot\hat{x}) \norm{L(v)x}^2 \,\dm(v),\end{equation}
where the norm above is the euclidean one.

Above $ L(v)\cdot\hat{x}$ denotes the projectivized element in  $ P(\C^k)$.

As  $\phi_L$ is stochastic we get that $\Pi_L(\hat{x}, P(\C^k)) = 1$. $\Pi_L(\hat{x}, S) $ describes the probability   of getting in the next step a state in  $S$, if the system is presently at the state $\hat{x}$.  \\

Remember that
$\tr(L(v)\pi_{\hat{x}} {L(v)}^\dagger) = \norm{L(v)x}^2 $, where $\pi_{\hat{x}} = \ket{x}\bra{x}$ and $x$ are representatives of norm $1$ in the class of $\hat{x}$.\\

This discrete-time process (described  by the kernel) taking values on $P(\C^k)$ is determined  by $\mu$ and $L$. If $\nu$ is a probability on the Borel $\sigma$-algebra $\mathcal{B}$ of $P(\C^k)$ define

\begin{align*}
		\nu\Pi_L(S) &= \int_{P(\C^k)} \Pi_L(\hat{x}, S) \,\, \dn(\hat{x}) \\
		&= \int_{P(\C^k) \times M_k} \textbf{1}_{S}(L(v)\cdot\hat{x}) \norm{L(v)x}^2 \,\, \dn(\hat{x})\dm(v).
\end{align*}

$\nu\Pi_L$ is a new probability on $ P(\C^k)$ and $\Pi_L$ is a Markov operator. The above definition of $\nu \to \nu \Pi_L$ is a simple generalization  of the one in $\cite{benoist2017invariant}$, where the authors take the $L$ considered here  as the identity transformation.\\

The map  $\nu \to \nu\,\Pi_L$ (acting on probabilities $\nu$) is called the Markov operator  obtained from  $\phi_L$  in the paper $\cite{lozinski2003quantum}$.
There the \textit{a priori} measure $\mu$ is a sum of Dirac probabilities. Here we consider a more general setting.

\begin{defi}
	We say that the probability $\nu$ over $P(\C^k)$  is invariant for  $\Pi_L$, if $\nu\Pi_L = \nu$.
\end{defi}

\medskip
The natural question is: does exist  such  invariant probability for $\Pi_L$ ?

\medskip

About the question of existence, we are going to prove that the kernel defined above is a continuous Markov operator (in the weak-star topology).
So, leaving the compact set of probabilities over $P(\C^k)$ invariant, by the Markov-Kakutani theorem there exists a fixed point, which means that there exists an invariant probability.
In order to do that we only need to find a linear operator $U: C_0(P(\C^k),\C ) \to C_0(P(\C^k),\C )$ such that $\langle Uf, \nu \rangle = \langle f,\nu \Pi_L \rangle$.
Here, $C_0(P(\C^k),\C )$ stands for continuous functions from $P(\C^k)$ to $\C$ with the $C_0$ norm which we denote by $\norm{\cdot}_{\infty}$.
When such $U$ exists we say that the Markov operator $\Pi_L$ is Feller.

According to Proposition 2.10 in $\cite{slomczynski2003dynamical}$ if such $U$ exists, then, $\Pi_{L}$ is continuous in weak-star topology and by Markov-Kakutani theorem, there is a fixed probability in $P(\C^k)$.

In Example $\ref{Mark1}$ we calculate the explicit expression of the invariant probability $\nu$.

\begin{theorem}\label{existence_inv_prob}
	Suppose that $L$ is such that $\int_{M_k} \norm{L(v)}^2\,\dm(v) < \infty$. Then, there exists at least one invariant probability $\nu$ for the Markov operator $\Pi_L$.
\end{theorem}

\begin{proof}
Define $U: C_0(P(\C^k),\C) \to C_0(P(\C^k),\C)$ by

\[		Uf(\hat x)=\int_{M_k} f(L(v)\cdot\hat{x}) \norm{L(v)x}^2 \,\dm(v).     \]

Notice that
\begin{align*}
 \langle Uf, \nu \rangle &= \int_{P(\C ^k)} Uf(\hat x)\, d\nu(\hat x)  \\
&= \int_{P(\C ^k) \times M_k} f(L(v)\cdot\hat{x}) \norm{L(v)x}^2,\dm(v) d\nu (\hat x)   \\
 &= \int_{P(\C ^k)} f(\hat x) \, d(\nu \Pi_L)(\hat x) =  \langle\nu \Pi_L\rangle.
\end{align*}

Therefore, $\langle Uf, \nu \rangle = \langle f,\nu \Pi_L \rangle$.

Then, we only need to prove that $Uf$ is a continuous function of $P(\C^k)$.

Consider a sequence $(\hat{x_n}) \in P(\C^k)$, such that, $\hat{x_n}\longrightarrow \hat x \in P(\C^k)$.
We are going to show that $Uf(\hat{x_n})\longrightarrow Uf(\hat x)$.
Define $F, F_n: M_k \to \mathbb{C}$ by
\[F_n(v)= f(L(v)\cdot \hat{x_n}) \norm{L(v)x_n}^2\]
and
\[F(v)= f(L(v)\cdot \hat{x}) \norm{L(v)x}^2\]

This way, $Uf(\hat{x_n})=\int F_n(v) \dm (v)$ and $Uf(\hat x)=\int F(v) \dm(v)$. Since the function $f$ and the norm are continuous, we have $F_n(v) \longrightarrow F(v)$, for all $v \in \mathcal{M}_k$.

Also,
$$ \abs{F_n(v)} = \abs{f(L(v)\cdot \hat{x_n})} \cdot \norm{L(v)x_n}^2 \le \norm{f}_{\infty} \ \tr (L(v) \ket{x_n} \bra{x_n} L(v)^{\dagger}) $$
$$ = \norm{f}_{\infty} \ \tr (\ket{x_n} \bra{x_n} L(v) L(v)^{\dagger}) \le \norm{f}_{\infty} \ \tr (L(v) L(v)^{\dagger})
= \norm{f}_{\infty} \norm{L(v)}^2. $$

As $\int \norm{L(v)}^2 \dm (v) < \infty$, we can apply Lebesgue Dominated Convergence Theorem to conclude that
$$ Uf(\hat{x_n})= \int F_n(v)\dm (v) \longrightarrow \int F(v) \dm (v) = Uf(\hat x).$$

So we have that $Uf$ is continuous and this is the end of the proof.
\end{proof}

\section{Process $\rho_n$, $n\in \mathbb{N}$, taking values on $\mathcal{D}_k$}\label{pro2}

For a fixed $\mu$ over $M_k$ and $L$ such that $\phi_L$ is irreducible and stochastic, one can naturally define a process $(\rho_n)$ on $\mathcal{D}_k = \{\rho\in M_k : \tr \rho = 1 \text{ and } \rho \ge 0\}$ which is called \emph{quantum trajectory} by T. Benoist, M. Fraas, Y. Pautrat, and C. Pellegrini in $\cite{benoist2017invariant}$. Given a $\rho_0$ initial state, we get
\[\rho_n = \frac{L(v)\rho_{n-1}{L(v)}^{\dagger}}{\tr(L(v)\rho_{n-1}{L(v)}^{\dagger})}\]
with probability
\[\tr(L(v)\rho_{n-1}{L(v)}^{\dagger})\dm(v),\,\,\,\,\,n \in \mathbb{N}.\]

This process has similarities with the previous one in  $P(\C^k)$ and we explore some relations between them. In this section, we follow closely  the notation of $\cite{benoist2017invariant}$.

We want to relate the invariant probabilities of the last section with the fixed point $\rho_{inv} = \rho_{inv}^L$  of $\phi_L$.

\medskip

First, denote $\Omega := M_k^\N$, and for $\omega={(\omega_i)}_{i\in\N}$, take $\pi_n(\omega) = (\omega_1,\ldots,\omega_n)$.
Recall that $\mathcal{M}$ is the Borel sigma-algebra on $M_k$.
For all, $n\in\N$, consider $\mathcal{O}_n$ the sigma algebra on $\Omega$ generated by the cylinder sets of size $n$, that is, $\mathcal{O}_n := \pi^{-1}_n(\mathcal{M}^{n})$.
We equip $\Omega$ with the smallest sigma algebra $\mathcal{O}$ which contains all $\mathcal{O}_n$, $n \in \mathbb{N}$.

Denote $\mathcal{J}_n := \mathcal{B}\times\mathcal{O}_n$ and $\mathcal{J} := \mathcal{B}\times\mathcal{O}$.
In this way, $(P(\mathbb{C}^k)\times\Omega, \mathcal{J})$ is an measurable space.
By abuse of language we  consider $V_i: \Omega \to M_k$ as a random variable  $V_i(\omega) = \omega_i$.
We also introduce another random variable $W_n := L(V_n)\cdots L(V_1)$, where $W_n(\omega) = L(\omega_n)\cdots L(\omega_1)$.\\

For a given a probability  $\nu$ on $P(\C^k)$, we define for  $S\in\mathcal{B}$ and $O_n\in\mathcal{O}_n$ another probability on $P(\C^{k}\times\Omega)$ by

\begin{equation}\label{Pnu}
\Pr_{\nu,n}(S\times O_n):=\int_{S \times O_n} \norm{W_n(\omega)x}^2\,\,\dn(\hat{x})\dm^{\otimes n}(\omega).
\end{equation}

\begin{remark}
	We can extend the above probability $\Pr_\nu$ over $\mathcal{B}\times\mathcal{O}$. We claim that $\Pr_{\nu,n}$, $n\in \mathbb{N}$,  is   a  consistent family over the cylinders of size $n$ (then, we can use the Caratheodory-Kolmogorov extension theorem).

		Indeed, note that $W_{n+1}(\omega)=L_{n+1}(\omega)W_{n}(\omega)$. Then
	\begin{align*}
		\Pr&_{\nu,n+1}(S\times O_n\times M_k) = \int_{S \times O_n\times M_k} \norm{W_{n+1}(\omega)x}^2\,\,\dn(\hat{x})\dm^{\otimes {n+1}}(\omega) \\
		&=\int_{S \times O_n \times M_k\hspace{-10px}} \tr\bigg(L(\omega_{n+1})W_{n}(\omega)\pi_{\hat x}W_n(\omega)^{\dagger} L(\omega_{n+1})^{\dagger}  \bigg)\,\,\dn(\hat{x})\dm^{\otimes {n+1}}(\omega)\\
		&=\int_{S \times O_n\hspace{-10px}} \tr\bigg(W_{n}(\omega)\pi_{\hat x}{W_n(\omega)}^{\dagger}{\int}_{M_k\hspace{-10px}} L(\omega_{n+1})^\dagger L(\omega_{n+1}) \dm(\omega_{n+1})\bigg)\,\dn(\hat{x})\dm^{\otimes n}(\omega)\\
		&=\int_{S \times O_n} \norm{W_n(\omega)x}^2 \,\,\dn(\hat{x})\dm^{\otimes n}(\omega)\\
		&=\Pr_{\nu,n}(S\times O_n).
	\end{align*}
\end{remark}

Since the set $\{W_n x =0\}$ leads to a null integrating term in (\ref{Pnu}), we have $\Pr_\nu(W_n x = 0) = 0$.  Therefore, we define the expression for each $n$ and then extend it. In this way $W_n(\omega) x \ne 0$. Remember that $W_n(\omega)\cdot \hat{x}$ is the representative of the class $W_n(\omega)x$, when $W_n(\omega)x \ne 0$.\\

\medskip

Denote $\E_\nu$ the expected value with respect to $\Pr_\nu$. Now observe that for a $\nu$ probability on $P(\C^k)$, if $\pi_{X_0}$ is an orthogonal projection on subspace generated by $X_0$ on $\C^k$, we have

\[ \rho_{\nu} := \E_\nu(\pi_{X_0}) = \int_{P(\C^k)} \pi_{x_0}\,\, \dn(x_0). \]

We call $\rho_{\nu}$ barycenter of $\nu$, and it is easy to see	that $\rho_\nu \in \mathcal{D}_k$.\\

Note that for each $\rho \in \mathcal{D}_k$, exists $(v_n)$ an orthonormal basis of eigenvectors with eigenvalues $a_i$ such that $\rho = \sum_i a_i \pi_{v_i}$. Therefore, exists $\nu = \sum a_i \delta_{v_i}$ such that $\rho_\nu = \rho$.\\

\medskip

We collect the above results in the next proposition (which was previously  stated as Proposition 2.1 in $\cite{benoist2017invariant}$ for the case $L=I$).

\medskip

\begin{proposition}\label{barycenter} If $\nu$ is invariant for $\Pi_L$, then
	$$\rho_\nu = \E_\nu(\pi_{\hat{X}_0}) = \E_\nu(\pi_{\hat{X}_1}) = \phi_L(\rho_\nu).$$

	Therefore,  for an irreducible $L$, every invariant measure $\nu$ for $\Pi_L$ has the same barycenter.
\end{proposition}

	We point out that in this way we can recover $\rho_{inv}$, the fixed point of $\phi_L$, by taking the barycenter of any invariant probability (the quantum channel $\phi_L$ admits only one fixed point). That is, for any invariant probability   $\nu$ for $\Pi_L$, we get that $\rho_{\nu}=\rho_{inv}$.

	Note that the previous process can be seen as $\rho_n:\Omega \to \mathcal{D}_k$, such that, $\rho_0(\hat{x},\omega) = \rho_\nu$ and, and $n \in \mathbb{N}$

$$\rho_n(\omega) =\frac{W_n(\omega)\rho_0 W_n(\omega)^{\dagger}}{\tr(W_n(\omega)\rho_0 W_n(\omega)^{\dagger})}.$$

Using an invariant $\rho$ we can define a Stationary Stochastic Process taking values on $M_k$.
That is, we will define a probability $\Pr$ over  $\Omega = (M_k)^\mathbb{N}$.

Take $O_n \in \mathcal{O}_n$ and define

$$\Pr^\rho (O_n) = \int_{O_n} \tr(W_n(\omega)\rho W_n(\omega)^{\dagger}) \,\,\dm^{\otimes n}(\omega).$$

The probability $\Pr$ on $\Omega$ defines a Stationary Stochastic Process.

\section{$\phi$-Erg and irreducible is Generic} \label{gen}
 \begin{defi}
      Given $L:M_k\to M_k$, $\mu$ on $M_k$ and $E$ subspace of $\C^k$,  we say that $E$ is $(L,\mu)$-invariant, if $L(v)E \subset E$, for all $v\in\supp\mu$.
    \end{defi}

    \begin{defi}
      Given $L: M_k\to M_k$, $\mu$ on $M_k$, we say that $L$ is $\phi$-Erg for $\mu$, if there exists an unique minimal non-trivial space $E$, such that, $E$ is $(L, \mu)$-invariant.
    \end{defi}

\medskip

In the case the space $E$ is equal to $\C^k$, as shown in Lemma $\ref{irred_E_Ck_equiv}$, we have $L$ irreducible for $\mu$ (or $\mu$-irreducible) in the sense of Definition $\ref{gil}$.

\medskip

Consider $\Bcal(M_k) = \{ L: M_k \to M_k \,\vert\, L \text{ is continuous and bounded} \}$ where $\norm{L} = \sup_{v\in M_k} \norm{L(v)}$. We write $\Bcal=\Bcal(M_k)$ when $k$ is implicit.

\begin{proposition}\label{perturbation_eigenvalues}
	Given $L \in \Bcal(M_k)$, $\mu$ over $M_k$, $v_1 \in $ supp $\mu$ and $\varepsilon > 0$, there exists $L_\varepsilon \in \Bcal(M_k)$ such that $\norm{L - L_\varepsilon} < \frac{\varepsilon}{2}$ and $L_\varepsilon(v_1)$ has $k$ distinct eigenvalues.
\end{proposition}

\begin{proof}
     Take $v_1 \in \supp\mu$. Denote by $J$ the Jordan canonical form for the complex matrix $L(v_1)$ and take $B$ such that $L(v_1) = B^{-1}JB$. Define $D_n = (d_{i,j})_{i,j} \in M_k$, where

    $$ d_{i,j} = \begin{cases}
      1 & \quad\text{if }i=n\,\,\text{and }j=n\\
      0 & \quad\text{otherwise.}
    \end{cases}$$

    Now, we look for each diagonal element of $J$. If the first, i.e., the element $(1,1)$ is zero, we sum $\frac{\delta}{4} D_1$. If the second element is not different from the first or is not different of zero, then, we sum $\frac{\delta}{2^i}D_2$, where $i > 2$ is chosen to satisfy both. We repeat this process until all the elements of diagonal are considered. After that, we get that all diagonal elements of $J+\sum_j \frac{\delta}{2^{i_j}} D_j$ are different and none is zero. Moreover, $\displaystyle\norm{\sum_j \frac{\delta}{2^{i_j}} D_j} \le \sum_j \frac{\delta}{2^{i_j}} \le \frac{\delta}{2}$.\\

    We define $\displaystyle D^\delta = \sum_j \frac{\delta}{2^{i_j}} D_j$ and $L_\varepsilon = L + B^{-1}D^\delta B$. Therefore, $\displaystyle\norm{ L_\varepsilon - L } = \norm{B^{-1}D^\delta B} \le \frac{\delta}{2}\norm{B^{-1}}\norm{B}$. Choosing $\delta < \frac{\varepsilon}{\norm{B^{-1}}\norm{B}}$ we get

    $$\norm{L_\varepsilon - L} < \frac{\varepsilon}{2}.$$

    Therefore, as $J+D^\delta$ has the same eigenvalues of $L_\varepsilon(v_1)$, we finished the proof.
\end{proof}

\begin{lemma}
  Consider  eigenvectors  $v_i \in \mathbb{C}^k, 1 \le i \le n$ of a linear transformation $A$ with respective eigenvalues $\lambda_i$, where $\lambda_i \neq \lambda_j$, for $i \neq j$.
  If a subspace $F \subseteq \mathbb{C}^k$ is invariant for $A$ and satisfies for some non-null constants $\alpha_1,\ldots,\alpha_n \in \mathbb{C}$
	\[\alpha_1v_1+\cdots+\alpha_nv_n \in F,\]

	then, $v_i \in F$ for all $1 \le i \le n$.
\end{lemma}

\textit{Proof.} We proceed by induction. Suppose $n=2$. Since $A(\alpha_1v_1+\alpha_2v_2) \in F$ and $ \lambda_1 (\alpha_1v_1+\alpha_2v_2) \in F$, we have
	$$ \lambda_1 (\alpha_1v_1+\alpha_2v_2) - A(\alpha_1v_1+\alpha_2v_2)$$
	$$ = \lambda_1 (\alpha_1v_1+\alpha_2v_2) - (\lambda_1 \alpha_1v_1+\lambda_2 \alpha_2v_2) $$
	$$= (\lambda_1 - \lambda_2)\alpha_2 v_2 \in F.$$

	Therefore, $v_1,v_2 \in F$. Now, assuming that the claim is true for every $n \le k$, we get
	$$ \lambda_{k+1} (\alpha_1v_1+\cdots+\alpha_{k+1}v_{k+1}) - A(\alpha_1v_1+\cdots+\alpha_{k+1}v_{k+1}) \in F.$$

	Which means $(\lambda_{k+1}-\lambda_1)\alpha_1v_1+\cdots+(\lambda_{k+1}-\lambda_k)\alpha_kv_k \in F$.
	From the hypothesis, this implies $v_1,\cdots,v_k \in F$. It follows that $v_{k+1} \in F$.

	\begin{theorem} \label{coi} Given $L \in \Bcal(M_k)$, $\mu$ over $M_k$ with $\#\supp\mu > 1$ and $\varepsilon > 0$, there exists $M_\delta \in \Bcal(M_k)$, such that, $\norm{L-M_\delta} < \varepsilon$ and $M_\delta$ is $\phi$-Erg and irreducible for $\mu$.
\end{theorem}

\medskip

\begin{proof}
  Given an $\varepsilon >0$, take $v_1\in \supp\mu$ such that $v_1 \neq 0$,  the respective $L_\varepsilon$ from Proposition $\ref{perturbation_eigenvalues}$ and moreover $\{x_1,\ldots,x_k\}$ such that  they are  a base of eigenvectors of $L_\varepsilon(v_1)$, with corresponding eigenvalues $\lambda_i$.
  If $L_\varepsilon$ is irreducible for $\mu$, we are done.
  Otherwise, there exists a decomposition in $E_1,\ldots, E_n$ minimal non-trivial subspaces that are invariant for all $L_\varepsilon(v)$, with $v$ in $\supp\mu$ and $k > \dim E_1 \ge \dim E_i$, for all $i$.

Remember that $E_i\cap E_j  = \{0\}$ and since all $E_i$ are invariant for $L_\varepsilon(v_1)$, they are generated by some of its eigenvectors.

Relabel $x_1,\ldots,x_k$ in such way that we get:

$E_1 = \langle x_1,\ldots,x_{d_1} \rangle$, $E_2 = \langle x_{d_1+1},\ldots,x_{d_2} \rangle , \ldots, E_n = \langle x_{d_{n-1}+1},\ldots,x_{d_n} \rangle $ and $K = \langle x_{d_n+1},\ldots,x_{k} \rangle$, with $\C^k = E_1 \oplus \cdots \oplus E_n \oplus K$, where $K$ is either $\{0\}$ or is not invariant for all $L_\varepsilon(v)$.

Now, define the linear transformation $A: \mathbb{C}^k \to \mathbb{C}^k$ by
$A(x_j)=x_{j+1}$. By abuse of notation, we assume that $x_{k+1}=x_1$.
Consider, for a $\delta >0$, the operator $M_\delta(v) = L_\varepsilon(v) + \frac{\delta \varphi(v)}{2 \norm{A}} A$, where $\varphi (v)= \frac{\norm{v-v_1}}{\norm{v}+\norm{v_1}}\le 1$.
Denote $c(v)=\frac{\delta \varphi(v)}{2\norm{A}} \ge 0$.
Note that $c(v) > 0$, for all $v \neq v_1$.
Notice that $M_\delta(v_1)=L_\varepsilon(v_1)$.
The idea here is to make an element $x_i$ move to all of the other subspaces, making it impossible to have an invariant and proper subspace for all $M_\delta(v)$. This combined with the proximity of the original $L$ will give us the result.

\medspace

\textit{Claim:} There exists a $\delta > 0$, such that the only non-trivial (and therefore minimal) subspace invariant for all $M_\delta(v)$, with $v \in supp \ \mu $, is $\C^{k}$.

Suppose $F \subseteq \C^{k}$ is such a subspace.
There  exists a non-trivial element $\alpha_1x_1+\cdots+\alpha_kx_k \in F \cap E_i$, for some constants $a_l \in \mathbb{C}^k$ and some $i$. This is so because if $K$ is $\{0\}$ or not invariant for $M_\delta(v_1)=L_\varepsilon(v_1)$, then $F \not\subset K$. Since not all $a_i$ can be zero, we have by the above lemma that some $x_j \in F$.

We take a matrix $v_2 \in supp \ \mu$, $v_2 \neq v_1$. Now,
$$M_\delta(v_2)x_j = L_\varepsilon(v_2)x_j+ c(v_2)Ax_j = L_\varepsilon(v_2)x_j+ c(v_2)x_{j+1}\in F.$$ As $E_i$ is invariant for $L_\varepsilon(v_2)$, we get that $$L_\varepsilon(v_2)x_j=\sum_{m=d_{i-1}+1}^{d_i} \alpha_m x_m.$$

Now, again, $F$ is invariant for $M_\delta(v_1)=L_\varepsilon(v_1)$, and then

$$ L_\varepsilon(v_1)M_\delta(v_2)(x_j)=L_\varepsilon(v_1)\left(\sum_{m=d_{i-1}+1}^{d_i}\alpha_m x_m + c(v_2)x_{j+1} \right)$$

$$ = \sum_{m=d_{i-1}+1}^{d_i}\lambda_m \alpha_m x_m + c(v_2)\lambda_{j+1} \ x_{j+1}  \in F.$$

Moving on, $L_\varepsilon(v_1)M_\delta(v_2)x_j - \lambda_{j+1} \cdot M_\delta(v_2)x_j \in F$. This means

$$ \sum_{m=d_{i-1}+1}^{d_i}(\lambda_m -\lambda_{j+1})\alpha_m x_m \in F.$$

By the lemma, $x_m \in F$, for all $m$ which are not $j+1$ and the corresponding $\alpha_m$ is not zero.
Now, suppose that $x_{j+1} \notin E_i$ (this excludes the possibility of $m=j+1$ above).
In this way, $\alpha_m x_m \in F$, for all $m \in \{d_{i-1}+1,\ldots,d_i\}$, with no exceptions. It follows that $\sum_m \alpha_m x_m \in F$ and

$$ M_\delta(v_2)x_j- \sum_{m=d_{i-1}+1}^{d_i} \alpha_m x_m \in F $$

$$ = \sum_{m=d_{i-1}+1}^{d_i} \alpha_m x_m + c(v_2)x_{j+1} - \sum_{m=d_{i-1}+1}^{d_i} \alpha_m x_m $$
$$ = c(v_2)x_{j+1} \in F.$$

As $c(v_2) \neq 0$, we get $x_{j+1} \in F$. Now suppose $x_{j+1} \in E_i$. Then

$$ M_\delta(v_2)x_j- \sum_{\substack{m=d_{i-1}+1 \\ m \neq j+1}}^{d_i} \alpha_m x_m \in F. $$

This means $c(v_2)x_{j+1}+\alpha_{j+1}x_{j+1} \in F$. If $c(v_2)+\alpha_{j+1}=0$ we get a problem. In order to fix this, we need that $\frac{\delta \varphi(v_2)}{2 \norm{A}} \neq -\alpha_{j+1} \iff \delta \neq \frac{-2\alpha_{j+1}\norm{A}}{\varphi(v_2)}$. But, note that  $\alpha_{j+1}$ does not depend on $\delta$. In fact, it appears only in the decomposition

$$L_\varepsilon(v_2)x_j=\sum_{m=d_{i-1}+1}^{d_i} \alpha_m x_m.$$

Since we can do this decomposition for all $j$, we only have to check that

$$ \delta \notin \left\{\frac{-2\alpha_{j+1}\norm{A}}{\varphi(v_2)}; 1\le j \le d_n\right\}.$$

Taking $\delta$ small enough, we accomplish this and also we get $\delta < \varepsilon$. Now, we get the claim in the same way: $x_{j+1} \in F$ and $F=\mathbb{C}^k$. So, for this $\delta$ we get that $M_\delta$ is irreducible. Finally,

$$ \norm{L-M_\delta} \le \norm{L-L_\varepsilon} + \norm{L_\varepsilon-M_\delta} < \varepsilon/2 + \norm{\frac{\delta \varphi(v)A}{2 \norm{A}}} < \varepsilon. $$
\end{proof}

    \begin{defi}
      For a fixed measure $\mu$ over $M_k$, define

      $$\Bcal_\mu(M_k) = \{ L\in \Bcal \,\vert\, \text{L irreducible for} \mu \},$$

      and

      $$\Bcal^\phi_\mu(M_k) = \{ L\in \Bcal \,\vert\, \text{L is $\phi$-Erg for $\mu$} \}.$$
    \end{defi}

    \begin{corollary}\label{phi_erg_dense_set}
      Given $\mu$ over $M_k$ with $\#\supp\mu > 1$, $\Bcal_\mu(M_k)$ is dense on $\Bcal(M_k)$.
    \end{corollary}
	\begin{proof}
		It follows from the above.
	\end{proof}

    \medskip

\begin{proposition}\label{irred_open}
	$\Bcal_\mu(M_k)$ is open for a fixed $\mu$ on $M_k$.
\end{proposition}

\begin{proof}
		We will prove that the complement of $\Bcal_{\mu}(M_k)$ is closed in $\Bcal(M_k)$. Let $L_n$ be a sequence outside $\Bcal_{\mu}(M_k)$ converging to some $L \in \Bcal(M_k)$. For each $n$, consider $E_n$ a non-trivial $(L_n,\mu)$-invariant  subspace and $P_n$ the projection on $E_n$.

	The $(L_n,\mu)$-invariance is equivalent to say that $L_n(v)P_n = P_nL_n(v)P_n$, for all $v\in\supp\mu$. Therefore, there is a subsequence such that $P_{n_i}\to P$, where $P$ is a projection. Rename $P_n\to P$. Furthermore, $L_n \to L$, thus $P_nL_n(v)P_n=L_n(v)P_n \to PL(v)P=L(v)P$, for all $v\in\supp\mu$. This implies that $E := \Im(P)$ is $(L,\mu)$-invariant for $L$. Of course, $E$ is not the trivial space because $\norm{P}\ge 1$. Moreover, we know that $\ker(P_n)$ is non-trivial for all $n$, once  $L_n$ is not $\mu$-irreducible. So, take $x_n\in \ker(P_n)$ with $\norm{x_n}=1$, and rename it in order to get a subsequence such that $x_n \to x$. Observe that $P_nx_n = 0$, for all $n$ and $P_nx_n \to Px$. This implies that $Px = 0$ and, of course, $\ker(P)$ is non-trivial. Hence, $E\ne \C^k$ and $L$ is not $\mu$-irreducible.

\end{proof}

    \begin{proposition}
       $\Bcal^\phi_\mu(M_k)$ is open for a fixed $\mu$ on $M_k$.
    \end{proposition}

	\begin{proof}
    Take $L_n \to L$ such that $L_n$ is not $\phi$-Erg. Therefore, there exists $E_{1,n} \oplus E_{2,n} \oplus E_{0, n} = \C^k$, with $E_{i,n}$ minimal $(L_n,\mu)$-invariant for $L_n$, where $i=1,2$ and $E_{0,n}$ is not necessarily $(L_n,\mu)$-invariant. Take $P_{i,n}$ the projection on $E_{i,n}$. Rename them in order to get a subsequence such that $P_{i,n}\to P_i$, for all $i=1,2,0$. By using the same argument as the one used in  Proposition $\ref{irred_open}$, we observe that $E_{i}=\Im(P_i)$ is $(L,\mu)$-invariant for $L$, for $i=1,2$. If $x\in E_1 \setminus \{0\}$ we know that $\lim_n \norm{P_{1,n}x - x} = \norm{P_1x - x} = 0$, so defining $x_n := P_{1,n}x \in E_{1,n}$, we get
    $x_n \to x$. As $0 = P_{2,n}x_n \to P_2x$, we know $x \in \ker P_2$ and therefore $x \notin E_2$. This argument shows that $E_1 \cap E_2 = \{0\}$, hence $L$ is not $\phi$-Erg because it admits two $(L, \mu)$-invariant subspaces.
	\end{proof}

    \begin{corollary}\label{corimp}
      Given $\mu$ over $M_k$ with $\#\supp\mu > 1$, $\Bcal^\phi_\mu(M_k)$ is open, dense and, therefore, generic.
    \end{corollary}

	\section{Some examples} \label{exam}

In this section, we present several examples. The main one is Example \ref{Mark1}  that considers a quantum channel which is a kind of  version of a Markov chain. We can show in expression \eqref{estea} that the entropy of this channel coincides with the entropy of the associated stationary Markov Process. This is a piece of clear evidence that  our definition is a natural extension of the classical  concept of entropy.  
In \cite{BKL} it is shown that the entropy of this channel is 
related to one of the Lyapunov exponents of the associated time evolution process which are described in sections \ref{pro1} and \ref{pro2}.

		\begin{example}

		Let $V_{2n}=c\cdot \left ( \begin{array}{cc} \frac{1}{2n} & 0 \\ 0 & 0 \end{array} \right)$
		and $V_{2n-1}=d\cdot \left( \begin{array}{cc} 0 & \frac{1}{2n-1} \\ 0 & 0 \end{array} \right)$,
		for all $n\ge 1$ (with constants $c$ and $d$ to be defined). Then,

		$$V_{2n}^{\dagger}V_{2n}=\frac{c^2}{(2n)^2} \cdot \left ( \begin{array}{cc} 1 & 0 \\ 0 & 0 \end{array} \right)
		\text{ and }\, V_{2n-1}^{\dagger}V_{2n-1}=\frac{d^2}{(2n-1)^2} \cdot \left( \begin{array}{cc} 0 & 0 \\ 0 & 1 \end{array} \right).  $$

		Setting $L=I$ (the identity map $v \mapsto v$) and $\mu = \sum_{n=1}^{\infty} \delta_{V_n}$, we have

		$$ \int_{M_k} L(v)^{\dagger}L(v) \dm (v) = \sum_{n=1}^{\infty} V_n^{\dagger}V_n$$
		$$ = c^2 \left ( \begin{array}{cc} 1 & 0 \\ 0 & 0 \end{array} \right) \sum_{n=1}^{\infty} \dfrac{1}{(2n)^2} + d^2 \left( \begin{array}{cc} 0 & 0 \\ 0 & 1 \end{array} \right) \sum_{n=1}^{\infty} \frac{1}{(2n-1)^2}.$$

		Choosing $$c = \left( \sum_{n=1}^{\infty} \frac{1}{(2n)^2} \right )^{-1/2} \text{ and } \, d = \left( \sum_{n=1}^{\infty} \frac{1}{(2n-1)^2}\right)^{-1/2},$$

		we get $\int L(v)^{\dagger}L(v) \dm (v) = \Id$. Now, notice that

		$$\int \norm{L(v)} \dm (v) = c \cdot \sum_{n=1}^{\infty} \frac{1}{2n}  + d \cdot \sum_{n=1}^{\infty} \frac{1}{2n-1} = \infty,$$
whereas $ \norm{L(v)} \le max \{c, d\} < \infty $, for all $v \in supp(\mu)$. Even when the last integral is not finite, the limitation on the norm above should produce an invariant probability for the kernel, according to Theorem $\ref{existence_inv_prob}$. To show this will be  our goal. Before that, we will compute the action of the quantum channel (in order to clear out what is the fixed density).

		 For a general density $\rho= \left( \begin{array}{cc} \rho_1 & \rho_2 \\ \rho_3 & \rho_4 \end{array} \right)$, we have

		 $$V_{2n}\rho V_{2n}^{\dagger}= \dfrac{c^2}{(2n)^2} \left( \begin{array}{cc} 1 & 0 \\ 0 & 0 \end{array} \right) \left( \begin{array}{cc} \rho_1 & \rho_2 \\ \rho_3 & \rho_4 \end{array} \right) \left( \begin{array}{cc} 1 & 0 \\ 0 & 0 \end{array} \right) = \dfrac{c^2}{(2n)^2} \left( \begin{array}{cc}\rho_1 & 0 \\ 0 & 0 \end{array} \right) , $$
		and
	\begin{align*}
	  V_{2n-1}\rho V_{2n}^{\dagger}&= \dfrac{d^2}{(2n-1)^2} \left( \begin{array}{cc} 0 & 1 \\ 0 & 0 \end{array} \right) \left( \begin{array}{cc} \rho_1 & \rho_2 \\ \rho_3 & \rho_4 \end{array} \right) \left( \begin{array}{cc} 0 & 0 \\ 1 & 0 \end{array} \right) \\
	&= \dfrac{d^2}{(2n-1)^2} \left( \begin{array}{cc} \rho_4 & 0 \\ 0 & 0 \end{array} \right).
	\end{align*}

		That is,

		$$ \phi_L(\rho)= \sum_{n=1}^{\infty} \left( \dfrac{c^2}{(2n)^2} \rho_1 +  \dfrac{d^2}{(2n-1)^2} \rho_4  \right) \left( \begin{array}{cc}1 & 0 \\ 0 & 0 \end{array} \right) = (\rho_1 + \rho_4) \left( \begin{array}{cc}1 & 0 \\ 0 & 0 \end{array} \right)$$
		$$ = \tr(\rho) \cdot \ket {e_1} \bra {e_1} = \ket {e_1} \bra {e_1}. $$

		\medskip

This $ \phi_L$ is not irreducible but it is an interesting example. It is a case where the invariant probability is unique as we will see soon.

		Clearly, the only fixed point for  $\phi_L$   is $\rho_{inv}=\ket {e_1} \bra {e_1}$.  What we should expect for invariant probabilities over $P(\C ^k)$? As the fixed point is itself a projection and the proposition $\ref{barycenter}$ says it is an average of projections around any invariant probability, the only option is a probability concentrated in $\hat e_1$, which is $\nu=\delta_{\hat e_1}$. Let's check that it is the case.

		For a general probability $\nu$ over $P(\C ^k)$ and a Borel set $B \subset P(\C ^k)$, we have

		$$ \nu \Pi_L(B)= \int_{M_k} \int_{P(\C ^k)} \textbf{1}_{B}(L(v)\cdot \hat x) \norm{L(v)x}_{HS}^2 \dm (v) d\nu(\hat x)$$
		$$ = \int_{P(\C ^k)} \sum_{n=1}^{\infty} \left[\textbf{1}_{B}(V_{2n}\cdot \hat x) \norm{V_{2n}x}_{HS}^2  + \textbf{1}_{B}(V_{2n-1}\cdot \hat x) \norm{V_{2n-1}x}_{HS}^2 \right] d\nu(\hat x).$$

		Notice that $V_{2n}\cdot \hat x = \hat e_1$ for $\hat x \neq \hat e_2$ and $V_{2n-1}\cdot \hat x = \hat e_1$ for $\hat x \neq \hat e_1$, whereas $V_{2n} e_1  = V_{2n-1} e_2 = 0$. Also, for a representative $x=(x_1,x_2)$ of norm 1, we got $(\ket x \bra x)_{ij}=x_i \overline{x_j}$. So,

	$$	\tr(V_{2n} \,\ket x \bra x \, V_{2n}^{\dagger})=\dfrac{c^2}{(2n)^2} \cdot (\ket x \bra x)_{11}= \dfrac{c^2}{(2n)^2} |x_1|^2 ,$$
and
	$$ \tr(V_{2n-1} \, \ket x \bra x \, V_{2n-1}^{\dagger})=\dfrac{d^2}{(2n-1)^2} \cdot (\ket x \bra x)_{22}= \dfrac{d^2}{(2n-1)^2} |x_2|^2.$$

		Then,

		$$ \nu \Pi_L(B)= \int_{P(\C ^k)} \sum_{n=1}^{\infty} \textbf{1}_{B}(\hat e_1) \left[ \dfrac{c^2}{(2n)^2} |x_1|^2  +\dfrac{d^2}{(2n-1)^2} |x_2|^2 \right] d\nu(\hat x) $$
		$$ = \int_{P(\C ^k)} \textbf{1}_{B}(\hat e_1) (|x_1|^2  + |x_2|^2 ) d\nu(\hat x) $$
		$$ = \int_{P(\C ^k)} \textbf{1}_{B}(\hat e_1) d\nu(\hat x) $$
		$$ = \textbf{1}_{B}(\hat e_1).$$

We conclude that if $\nu \Pi_L=\nu$, then $\nu=\delta_{\hat e_1}$. We also get a bonus: the invariant probability is unique.

		\end{example}

		To illustrate Proposition $\ref{barycenter}$ (under the irreducible condition)  we write down the following example.

		\begin{example}

The next example is  somehow related to Example $\ref{Mark1}$. Let's define

		$$ V_1=\left ( \begin{array}{cc} 1 & 0 \\ 0 & 0 \end{array} \right ) \text{ and } V_2=\left ( \begin{array}{cc}  0 & 1 \\ 0 & 0 \end{array} \right ). $$

		These two matrices  generate the same elements which we will consider in Example $\ref{Mark1}$, since for $\mu=\delta_{V_1}+\delta_{V_2}$,

		$$ \phi_I(\rho)=V_1 \rho V_1 ^{\dagger} + V_2 \rho V_2 ^{\dagger} = \ket{e_1} \bra{e_1}.$$

		Also, we get that $\phi_I$ is not irreducible. Wanting to fix this issue, we introduce

$$ V_3=\left ( \begin{array}{cc} 0 & 0 \\ 1 & 0 \end{array} \right ) \text{ and } V_4=\left ( \begin{array}{cc}  0 & 0 \\ 0 & 1 \end{array} \right ). $$

		Notice that these two matrices generates another channel $\psi$ that maps every density $\rho$ into $\ket{e_2} \bra{e_2}$. So, it is also not irreducible.
		Now, redefining $\mu = \frac{1}{2} \sum_{i=1}^4 \delta_{V_i}$, we get that

		$$\phi_I(\rho)= \frac{1}{2} \sum_{i=1}^4 V_i \rho V_i ^{\dagger} = \frac{1}{2} (\ket{e_1}\bra{e_1}+\ket{e_2}\bra{e_2}) = \frac{1}{2} \Id. $$

In this case, $\mu$ is a measure and not a probability.

		We compute the products

		$$ V_1^{\dagger}V_1=V_1, \, V_2^{\dagger}V_2=V_4,$$
and
		$$ V_3^{\dagger}V_3=V_1 \text{ and } \, V_4^{\dagger}V_4=V_4.$$

		In this way,
		$$ \phi_I^{*}(\Id)= \frac{1}{2} \sum_{i=1}^4 V_i ^{\dagger} V_i = V_1 + V_4 = \Id,$$
and $\phi_I$ is stochastic. As $\Id>0$, we get that $(I+\phi)(\rho)=\rho+\phi(\rho)=\rho+\Id>0$, and so $\phi$ is irreducible. Clearly, $\rho_{inv}=\frac{1}{2}\, \Id$.

		Now, for a general $\nu$ over $P(\C ^k)$ and a Borel set $B \subset P(\C ^k)$, we get

	$$ \nu \Pi_I(B)= \int_{P(\C ^k)} \int_{M_k} \textbf{1}_{B}(L(v)\cdot \hat x) \norm{L(v)x}_{HS}^2 \, \dm (v) d\nu(\hat x)$$
	$$ =  \int_{P(\C ^k)} \sum_{i=1}^4 \frac{1}{2} \textbf{1}_{B}(V_i \cdot \hat x) \norm{V_i x}_{HS}^2 \, d\nu(\hat x).$$

Remember that

$$  V_1 \, \ket x \bra x \, V_1^{\dagger}= \left(\begin{array}{cc} |x_1|^2 & 0 \\ 0 & 0
		\end{array}\right), \ \
		V_2 \, \ket x \bra x \, V_2^{\dagger}= \left(\begin{array}{cc} |x_2|^2 & 0 \\ 0 & 0
		\end{array}\right),$$
$$ V_3 \, \ket x \bra x \, V_3^{\dagger}= \left(\begin{array}{cc} 0 & 0 \\ 0 & |x_1|^2
		\end{array}\right) \ \text{ and } \
		V_4 \, \ket x \bra x \, V_4^{\dagger}= \left(\begin{array}{cc} 0 & 0 \\ 0 & |x_2|^2
		\end{array}\right).$$

So,

$$ \nu \Pi_I(B)=  \frac{1}{2} \int_{P(\C ^k)}  [\textbf{1}_{B}(V_1 \cdot \hat x) + \textbf{1}_{B}(V_3 \cdot \hat x)] |x_1|^2 + [\textbf{1}_{B}(V_2 \cdot \hat x) + \textbf{1}_{B}(V_4 \cdot \hat x)] |x_2|^2 \,d\nu(\hat x)$$
$$=\frac{1}{2} \int_{P(\C ^k)}  [\textbf{1}_{B}(\hat e_1) + \textbf{1}_{B}(\hat e_2)] |x_1|^2 + [\textbf{1}_{B}(\hat e_1) + \textbf{1}_{B}(\hat e_2)] |x_2|^2 \, d\nu(\hat x)$$
$$=\frac{1}{2} \int_{P(\C ^k)}  \textbf{1}_{B}(\hat e_1) + \textbf{1}_{B}(\hat e_2) \, d\nu(\hat x)$$
$$= \frac{1}{2} \textbf{1}_{B}(\hat e_1) + \textbf{1}_{B}(\hat e_2)$$
$$= \frac{1}{2} \, \delta_{\hat e_1}(B) + \frac{1}{2} \, \delta_{\hat e_2}(B).$$

We conclude that if $\nu= \nu \Pi_I$, then $\nu = \frac{1}{2} \, \delta_{\hat e_1} + \frac{1}{2} \, \delta_{\hat e_2}$. Note that (see the concept of barycenter in Section $\ref{pro2}$)

$$\int_{P(\C ^k)} \pi_{x}\, d\nu(\hat x) = \frac{1}{2} \,\pi_{e_1}+\frac{1}{2}\, \pi_{e_2}=\frac{1}{2} \Id = \rho_{inv}.$$

		\end{example}

		\begin{example}[$L$ is a  $C^*$-automorphism]
			Suppose that  $\mu$ over $M_k$ satisfies the below conditions:
			\begin{itemize}
				\item $\displaystyle \int_{M_k} v^\dagger v \, \dm(v) = \Id; $ and
				\item $\displaystyle \int_{M_k} \norm{v}^2 d\mu(v) < \infty$, where $\norm{\cdot}$ is the Hilbert-Schmidt norm.
			\end{itemize}

			Take an unitary matrix $U\in M_k$ and define $L(v) = UvU^\dagger$. Note that $\norm{UvU^\dagger}^2 = \tr(UvU^\dagger) = \tr(v) = \norm{v}^2$. Moreover,

			\begin{align*}
			\int_{M_k} L(v)^\dagger L(v) \, \dm(v) &= \int_{M_k} Uv^\dagger U^\dagger UvU^\dagger \, \dm(v)\\
			&= U \int_{M_k} v^\dagger v \, \dm(v) U^\dagger \\
			&= \Id.
			\end{align*}
		\end{example}

		\begin{remark}
			The operators of the form $L(v)=UvU^\dagger$, for an unitary $U$, are the  $C^*$-automorphisms of $M_k$ (see Section 1.4 in $\cite{arveson1998invitation}$).
		\end{remark}

In the next example, we adapt the reasoning of an Example $4$ in $\cite{baraviera2010thermodynamic}$ to the present setting.

We will show that for a certain
$\mu$
and
$L$
(and, quantum channel) the value we get here for the entropy is equal to the classical entropy of a Markov Chain (when the state space is finite).

	\begin{example}[The Markov model in quantum information]\label{Mark1}
		Suppose that $P=\left(\begin{array}{cc}
		p_{00} & p_{01} \\ p_{10} & p_{11}
		\end{array}\right)$ is a irreducible (in the classical sense for a Markov chain)  column stochastic matrix. Define $\mu$ over  $M_2$ by
		$$\mu = \sum_{i=1}^{4} \delta_{V_i},$$
where the matrices $V_i$ are

		$$ V_1= \left(\begin{array}{cc} \sqrt{p_{00}} & 0 \\ 0 & 0
		\end{array}\right), \
		V_2= \left(\begin{array}{cc} 0 & \sqrt{p_{01}} \\ 0 & 0
		\end{array}\right),$$

		$$ V_3= \left(\begin{array}{cc} 0 & 0 \\ \sqrt{p_{10}} & 0

		\end{array}\right) \text{ and }
		V_4= \left(\begin{array}{cc} 0 & 0 \\ 0 & \sqrt{p_{11}}
		\end{array}\right).		$$

		We take $L=I$ and $\phi_I=\phi_L$, in order to get the quantum channel

		$$\phi(\rho) = \sum_{1}^{4} V_i\rho V_i^{\dagger},$$

		whose dual is

		$$\phi^{*}(\rho)=\sum_{1}^{4} V_i^{\dagger} \rho V_i.$$

		Note that
		$$ V_1^{\dagger}V_1= \left(\begin{array}{cc} p_{00} & 0 \\ 0 & 0
		\end{array}\right), \ \
		V_2^{\dagger}V_2= \left(\begin{array}{cc} 0 & 0 \\ 0 & p_{01}
		\end{array}\right)$$

		\begin{equation}\label{equacoes_VTV}
		V_3^{\dagger}V_3= \left(\begin{array}{cc} p_{10} & 0 \\ 0 & 0
		\end{array}\right) \ \,\text{and}\,\, \
		V_4^{\dagger}V_4= \left(\begin{array}{cc} 0 & 0 \\ 0 & p_{11}
		\end{array}\right),
		\end{equation}

		that is,

		$$ \phi^{*}(\Id_2)=\left(\begin{array}{cc} p_{00}+ p_{10} & 0 \\ 0 & p_{01} + p_{11}
		\end{array}\right) = \Id_2 $$

		The channel $\phi$ is stochastic. We claim that the channel is irreducible (later we will exhibit the associated invariant density operator  $\rho$). Consider first the positive operator

		$$ \rho = \left(\begin{array}{cc} \rho_1 & \rho_2 \\ \rho_3 & \rho_4
		\end{array}\right)$$
	where $\rho_1,\rho_2 \in \mathbb{R} \text{ and } \ \rho_3=\overline{\rho_2}$ (in order to get that $\rho \ge 0$)

	The $V_i \rho V_i^{\dagger}$ are given by:

		$$  \rho^1 := V_1\rho V_1^{\dagger}= \left(\begin{array}{cc} p_{00}\rho_1 & 0 \\ 0 & 0
		\end{array}\right), \ \
		\rho^2 := V_2\rho V_2^{\dagger}= \left(\begin{array}{cc} p_{01}\rho_4 & 0 \\ 0 & 0
		\end{array}\right)$$

		\begin{equation}\label{equacoes_VrhoVT}
		 \rho^3 := V_3\rho V_3^{\dagger}= \left(\begin{array}{cc} 0 & 0 \\ 0 & p_{10}\rho_1
		\end{array}\right) \ \,\text{and}\,\, \
		\rho^4:= V_4\rho V_4^{\dagger}= \left(\begin{array}{cc} 0 & 0 \\ 0 & p_{11}\rho_4
		\end{array}\right)	\ \ \ \ \
		\end{equation}

		It follows that

		$$\phi(\rho)=\left(\begin{array}{cc} p_{00}\rho_1 + p_{01}\rho_4 & 0 \\ 0 & p_{10}\rho_1 + p_{11}\rho_4
		\end{array}\right).
		$$

		In the diagonal one can find the classical action  on vectors  of the Markov Chain described by $P$.

		In the same way for $v=(v_1,v_2) \in \mathbb{C}^2$, we get

$$\braket{v}{\phi(\rho)v}= (p_{00}\rho_1 + p_{01}\rho_4)|v_1|^2+ (p_{10}\rho_1 + p_{11}\rho_4)|v_2|^2 \ge 0.$$

		Moreover, the equality only happens when
		$$ p_{00}\rho_1 + p_{01}\rho_4=p_{10}\rho_1 + p_{11}\rho_4=0.$$

From this we get $\rho_1=\rho_4=0$, because $p_{ij}\geq 0$.

In this case, we get $\rho = 0$.

This means that , $\rho \neq 0, \rho \ge 0 \Rightarrow \phi(\rho)>0$, and, finally, we get that  $\phi$ is positive improving. From this, it follows that  $\phi$ is irreducible.

Now, we will look for the invariant density matrix.  Assuming $\rho_1+\rho_4=1$, we observe that  $\phi(\rho)=\rho \Rightarrow \rho_2=\rho_3=0$, and
		\begin{equation}\label{sistema}
		\left\{ \begin{array}{cc}
		\rho_1= p_{00}\rho_1 + p_{01}\rho_4 \\
		\rho_4=p_{10}\rho_1 + p_{11}\rho_4.

		\end{array}  \right.	\end{equation}

		We get $$(1-p_{00})\rho_1=p_{01}\rho_4=p_{01}(1-\rho_1)=p_{01}-p_{01}\rho_1 $$
		$$\Rightarrow (1-p_{00}+p_{01})\rho_1 = p_{01}.$$

		As $P$ is irreducible, it follows that $0<p_{ij}<1$ e
		$1-p_{00}+p_{01}>0$. That is,

		$$\rho_1=\dfrac{p_{01}}{1-p_{00}+p_{01}} \ \text{ and } \ \rho_4=\dfrac{1-p_{00}}{1-p_{00}+p_{01}}. $$

		An invariant density matrix is

		$$ \rho = \left(\begin{array}{cc} \dfrac{p_{01}}{1-p_{00}+p_{01}}   & 0 \\ 0 & \dfrac{1-p_{00}}{1-p_{00}+p_{01}}
		\end{array}\right).$$

		Note that $\pi=(\rho_1, \rho_4)\in \mathbb{R}^2$ is the  vector of probability which is invariant for the stochastic matrix $P$ (see  ($\ref{sistema}$)).

		Now, we will estimate the entropy of the quantum channel $\phi$.
			Using ($\ref{equacoes_VrhoVT}$) in the expression $\tr(V_jV_i\rho V_i^{\dagger}V_j^{\dagger})$ we get

			$$\left\{ \begin{array}{cccc}
			\tr(V_1\rho^i V_1^{\dagger})=p_{00}(\rho^i)_1 \\
			\tr(V_2\rho^i V_2^{\dagger})=p_{01}(\rho^i)_4 \\
			\tr(V_3\rho^i V_3^{\dagger})=p_{10}(\rho^i)_1 \\
			\tr(V_4\rho^i V_4^{\dagger})=p_{11}(\rho^i)_4

			\end{array} \right.$$

			For example,  $$\tr(V_3 V_1 \rho V_1^{\dagger} V_3^{\dagger})=\tr(V_3\rho^1 V_3^{\dagger})=p_{10}(\rho^1)_1 = p_{10}p_{00}\rho_1.$$

			From this we get the table.
			\begin{table}[h]
				\centering
				\vspace{0.5cm}
				\begin{tabular}{|c|cc|c|c|c|}
					\hline
					$ \tr(V_jV_i\rho V_i^{\dagger} V_j^{\dagger})$  &  $i$ & 1 & 2 & 3 & 4  \\ 
					\hline
					$j$ &  &  &  &  & \\
					1 & & $p_{00}^2\rho_1$	 & $p_{00}p_{01}\rho_4$ & 0 &0 \\
					2 & &		0	 &		0	 &		$p_{01}p_{10}\rho_1	$	 &	$p_{01}p_{11}\rho_4$		 \\
					3 &	&	$p_{00}p_{10}\rho_1	$	 &	$p_{10}p_{01}\rho_4 $  &	0		 &	0		 \\
					4 &	&		0	 &		0	 &		$p_{11}p_{10}\rho_1$	 &	$p_{11}^2\rho_4$	 \\
					\hline
					$\tr(V_i \rho V_i^{\dagger})$ &  &	$p_{00}\rho_1$	 &	$p_{01}\rho_4$	 &	$p_{10}\rho_1$	 &  	$p_{11}\rho_4$\\
					\hline
				\end{tabular}
			\end{table}

			The entropy we defined in the text is given by

			$$ h_{\mu}(L)=-\int_{M_k \times M_k} \tr(L(v)\rho L(v)^{\dagger})P(v,w)\log(P(v,w))d\mu(v)d\mu(w), $$

			where $P(v,w)= \frac{\tr(L(w)L(v)\rho L(v)^{\dagger}L(w)^{\dagger})}{\tr(L(v)\rho L(v)^{\dagger})}$.

We assumed before that $L=I$ and $\mu=\sum_i \delta_{V_i}$. Then, we finally get,

			$$ h_{\mu}(I)= -\sum_{i=1}^4 \sum_{j=1}^4 \tr(V_jV_i\rho V_i^{\dagger} V_j^{\dagger})\cdot \log \left(\dfrac{\tr(V_jV_i\rho V_i^{\dagger} V_j^{\dagger})}{\tr(V_i\rho V_i^{\dagger})}\right) $$

			$$ = -\big[p_{00}^2 \rho_1 \\log(p_{00})+p_{00}p_{10}\rho_1 \log(p_{10})+ p_{00}p_{01}\rho_4 \log(p_{00})+ p_{10}p_{01}\rho_4\log(p_{10}) $$
			$$ + p_{01}p_{10}\rho_1 \log(p_{01})+p_{11}p_{10}\rho_1 \log(p_{11})+p_{01}p_{11}\rho_4 \log(p_{01})+p_{11}^2\rho_4 \log(p_{11})\big] $$

			$$ = -\big[p_{00}\log(p_{00})(p_{00}\rho_1+p_{01}\rho_4)+p_{10} \log(p_{10})(p_{00}\rho_1+p_{01}\rho_4) $$
			$$+ p_{01} \log(p_{01})(p_{10}\rho_1+p_{11}\rho_4)+p_{11}\log(p_{11})(p_{10}\rho_1+p_{11}\rho_4) \big]$$

			$$ = -p_{00}\log(p_{00})\rho_1 - p_{10} \log(p_{10})\rho_1 -  p_{01} \log(p_{01})\rho_4 - p_{11}\log(p_{11})\rho_4$$
			$$ = - p_{00}\log(p_{00})\pi_0 - p_{10} \log(p_{10})\pi_0 -  p_{01} \log(p_{01})\pi_1 - p_{11}\log(p_{11})\pi_1= $$
$$ -\,\sum_{i,j=0}^1 \pi_j p_{ij} \log (p_{ij}).$$

Therefore, 
		\begin{equation} \label{estea}   h_{\mu}(I)= -\,\sum_{i,j=0}^1 \pi_j p_{ij} \log (p_{ij}).
		\end{equation}

			The last expression is the value of the classical Shannon-Kolmogorov entropy of the stationary Markov  Process  associated to the line stochastic matrix $P=(p_{ij})_{i,j=0,1}$ (see $\cite{Sp}$ and $\cite{PY}$).

The entropy is positive because the a priori $\mu$ is a measure (of mass equal to $4$) and not a probability.

			Now, let's look at the kernel $\Pi_L$ and find an invariant probability. For a given probability $\nu$ in $P(\C ^k)$ and a Borel set $B \subset P(\C ^k)$, we have

$$ \nu \Pi_L(B)= \int_{P(\C ^k)} \int_{M_k} \textbf{1}_{B}(L(v)\cdot \hat x) \norm{L(v)x}_{HS}^2  \, \dm (v) d\nu(\hat x),$$
which means
$$ \nu \Pi_L(B)= \int_{P(\C ^k)} \sum_{i=1}^4 \textbf{1}_{B}(V_i \cdot \hat x) \norm{V_i \, x}_{HS}^2 \, d\nu(\hat x).$$

Note that
$$ V_1\cdot \hat x=\hat e_1, \text{ if } \hat x \neq \hat e_2; \ \ V_2\cdot \hat x=\hat e_1, \text{ if } \hat x \neq \hat e_1;$$
$$ V_3\cdot \hat x=\hat e_2, \text{ if } \hat x \neq \hat e_2; \ \  V_4\cdot \hat x=\hat e_2, \text{ if } \hat x \neq \hat e_1 $$
$$ \text{ and } \, V_1(e_2)=V_2(e_1)=V_3(e_2)=V_4(e_1)=0.$$

It follows that

$$ \nu \Pi_L(B)= \int_{P(\C ^k)} \textbf{1}_{B}(\hat e_1) \, [\norm{V_1 x} + \norm{V_2 x}] + \textbf{1}_{B}(\hat e_2) \,  [\norm{V_3 x} + \norm{V_4 x}]   \, d\nu(\hat x).$$

			Now, we compute

		$$ \tr(V_1 \, \ket x \bra x \, V_1^{\dagger})= p_{00} \, |x_1|^2,$$
		$$ \tr(V_2 \, \ket x \bra x \, V_2^{\dagger})= p_{01} \, |x_2|^2,$$
		$$ \tr(V_3 \, \ket x \bra x \, V_3^{\dagger})= p_{10} \, |x_1|^2 $$
		$$ and \ \, \tr(V_4 \, \ket x \bra x \, V_4^{\dagger})= p_{11} \, |x_2|^2.$$

			In this way, we get

$$ \nu \Pi_L(B)= \int_{P(\C ^k)} \textbf{1}_{B}(\hat e_1) \, (p_{00}\, |x_1|^2 + p_{01}\, |x_2|^2 ) + \textbf{1}_{B}(\hat e_2) \,  (p_{10}\, |x_1|^2 + p_{11}\, |x_2|^2)   \, d\nu(\hat x).$$

From the last expression, we conclude that $\nu \Pi_L$ has support in the set $\{ \hat e_1, \hat e_2\}$.

In this way, if $\nu=\nu \Pi_L$, then it has to be equal to $\alpha \cdot \delta_{\hat e_1} + \beta \cdot \delta_{\hat e_2}$, with constants $\alpha, \beta \ge 0$, such that, $\alpha + \beta =1$. As we know the expression for $\rho_{inv}$, we can go further:

$$ \rho_{inv}=\int_{P(\C ^k)} \pi_{x} \, d\nu(\hat x) = \alpha \cdot \pi_{e_1}+\beta \cdot \pi_{e_2}.$$

As $$\rho_{inv}= \left ( \begin{array}{cc} \dfrac{p_{01}}{1-p_{00}+p_{01}} & 0 \\ 0 & \dfrac{1-p_{00}}{1-p_{00}+p_{01}} \end{array}	 \right ),$$
we get that $\alpha=\dfrac{p_{01}}{1-p_{00}+p_{01}}$ and $\beta=\dfrac{1-p_{00}}{1-p_{00}+p_{01}}$.
\medskip

In order to finish our example, we write down the invariant probability
$$ \nu = \dfrac{p_{01}}{1-p_{00}+p_{01}} \cdot \delta_{\hat e_1}+ \dfrac{1-p_{00}}{1-p_{00}+p_{01}}\cdot \delta_{\hat e_2} = \pi_1\, \delta_{\hat e_1} + \pi_2 \, \delta_{\hat e_2} ,$$
and we point out that the two constants are no more no less then the entries of the invariant probability vector $\pi=(\pi_1, \pi_2)$ for the Markov chain with transitions $P=(p_{ij})_{i,j=1,2}$.

\medskip

 In this way, the concept of entropy we considered before in Section $\ref{ent}$ is a natural generalization
 of the classical Kolmogorov-Shannon entropy and the process $X_n, n \in \mathbb{N},$ of Section $\ref{pro1}$ is a natural generalization of the classical Markov Chain process.

\medskip

\end{example}

\medskip

\begin{example}\label{gre}
Consider  a measure $\mu$  with support on the set
$$\{ \left ( \begin{array}{cc} x & -y \\ y& x \end{array}	 \right ) |\, x,y \in \mathbb{R}\}\, \subset M_2,$$
such that has density  $ f(x,y)= \frac{1}{4 \pi} \, e^{- \,\frac{(x^2 \,+\, y^2)}{2}}$ (see also (9) in $\cite{Wo}$)

Taking $L=I$ we get that $\rho_0= \left ( \begin{array}{cc} 1/2 & 0 \\ 0& 1/2 \end{array}	 \right )$ satisfies
$\phi_I (\rho_0)=\rho_0$.

Indeed the channel is given by
$$ \rho= \left ( \begin{array}{cc} a & b \\ c& d \end{array}	 \right ) \, \to \, \phi_I(\rho)=$$
$$\int \int  \, \left ( \begin{array}{cc} x & - y \\ y& x \end{array}	 \right ) \left ( \begin{array}{cc} a & b \\ c& d \end{array}	 \right )\left ( \begin{array}{cc}  x& y \\ -y & x \end{array}	 \right )  \frac{1}{4 \pi} \, e^{- \,\frac{(x^2 \,+\, y^2)}{2}}\,dx\,dy= \left ( \begin{array}{cc} 1/2 & \frac{b-c}{2} \\\frac{c-b}{2}& 1/2 \end{array}	 \right ). $$

\medskip
Notice that although $\left (\begin{array}{cc} 1/2 & b \\\ -b & 1/2 \end{array}	 \right )$ is a fixed point of $\phi_I$, it is not a density unless $b=0$. Thus, $\rho_0$ is the only eigendensity.

Given a probability  $\nu$ on $P(\C^k) $  the expression for the kernel is

		$$
			 \nu\Pi_L(S) = \int_{P(\C^k)} \Pi_L(\hat{w}, S) \,\, \dn(\hat{w})
			 	=$$
 $$\int_{P(\C^k) \times M_k} \textbf{1}_{S}(L(v)\cdot\hat{w}) \norm{L(v)w}^2 \,\, \dn(\hat{w})\dm(v)=
		$$
		$$\int_{P(\C^k) \times M_k} \textbf{1}_{S} \widehat{\left ( \begin{array}{cc} v_1\, w_1 - v_2\, w_2 \\ v_2\, w_1 + v_1\, w_2 \end{array}	 \right )}\,\, (v_1^2 + v_2^2)  \,\frac{1}{4 \pi} \, e^{- \,\frac{(v_1^2 \,+\, v_2^2)}{2}}\, dv_1\, d v_2 \,\dn(\hat{w}).
		$$
\medskip

Now, we will estimate the entropy (which will be negative).

\medskip

Using the fixed density operator $\rho_0=\matrixx{\frac{1}{2} & 0\\0 & \frac{1}{2}}$ we get (according to Section $\ref{ent}$)

$$P(v,w) = \frac{\tr(wv \rho_0 v^\dagger w^\dagger)}{\tr(v \rho_0 v^\dagger)}.$$

We denote

$$ w = \matrixx{w_1 & -w_2\\ w_2 & w_1}\text{ and } v = \matrixx{v_1 & -v_2\\ v_2 & v_1},$$

and we get

\begin{align*}
	\tr(v \rho_0 v^\dagger) &=\frac{1}{2} \tr\left(\matrixx{v_1 & -v_2\\ v_2 & v_1} \matrixx{v_1 & v_2\\ -v_2 & v_1}    \right) \\
	&=\frac{1}{2}\tr \matrixx{v_1^2 + v_2^2 & 0\\0&v_1^2+v_2^2} = v_1^2 + v_2^2
\end{align*}

and

\begin{align*}
	\tr(wv \rho_0 v^\dagger w^\dagger) &= \frac{1}{2} \tr\left(w\matrixx{v_1 & -v_2\\ v_2 & v_1} \matrixx{v_1 & v_2\\ -v_2 & v_1}  w^\dagger  \right) \\
	&= (v_1^2 + v_2^2)\tr(w w^\dagger) = (v_1^2 + v_2^2)(w_1^2+w_2^2).
\end{align*}

Thus, we get the following expression for the entropy (remember that $\int_0^\infty x^3 e^{-\frac{x^2}{2}} dx = 2$):

\begin{align*}
	h_\mu(L) &= - \frac{1}{16\pi^2}\int (v_1^2 + v_2^2)(w_1^2 + w_2^2) \log(w_1^2 + w_2^2) e^{-\frac{v_1^2+v_2^2}{2}} e^{-\frac{w_1^2+w_2^2}{2}} dv_1 dv_2 dw_1 dw_2\\
	&= - \frac{1}{4}\int_0^{\infty}\int_0^{\infty} r_v^3 r_w^3 \log(r_w^2) e^{-\frac{r_v^2}{2}}e^{-\frac{r_w^2}{2}} dr_v dr_w \\
	&= -\frac{1}{4}\int_0^{\infty} \left[\int_0^{\infty} r_v^3 e^{-\frac{r_v^2}{2}} dr_v \right] r_w^3\log(r_w^2)e^{-\frac{r_w^2}{2}} dr_w\\
	&= -\frac{1}{2}\int_0^{\infty} r_w^3\log(r_w^2)e^{-\frac{r_w^2}{2}} dr_w\\
	&= - \int_{0}^{\infty} r_w^3 \log(r_w) e^{-\frac{r_w^2}{2}} dr_w\\
	&\approx -1.11593
\end{align*}

We used polar coordinates above.

\end{example}

\section{Conclusion and relations with other works} \label{Con}

We introduce a concept of entropy and pressure (definitions depending on an {\it a priori} probability $\mu$).
For  a given  $H: M_k \to M_k$ (which plays the role of an Hamiltonian, or a Liouvillian)  we define a version of the Ruelle operator
$\phi_H:M_k \to M_k$, via the expression:
\[\rho \,\to\,\phi_H(\rho) = \int_{M_k} H(v) \rho {H(v)}^\dagger \, \dm(v).\]

\medskip

 After that, we presented a type of  Ruelle Theorem: a variational principle  of pressure related to an eigenvalue problem for the  Ruelle operator (see Theorem \ref{peqlog}).  The entropy and the Ruelle operator are linked via the {\it a priori} probability in a natural and fundamental way.

The definition of entropy considered here is not based on the point of view of dynamical partitions. It is a kind of generalization of  Rokhlin Formula which says the entropy of an $\sigma$-invariant probability $\nu$ is $H(\nu)=-  \int \log J d \nu$, where $J$ is the Jacobian (a dynamical version of Radon-Nikodym derivative). Note that this entropy is not relative but absolute. Results in \cite{LMMS}   - for the classical (not quantum) Thermodynamic Formalism theory -  include the case where the alphabet $M$  (a compact metric space) is uncountable. We did not use the results of \cite{LMMS} we just mentioned it to say  that we followed similar reasoning.

A common procedure in Statistical Mechanics (for the one-dimensional lattice
$M^\mathbb{N}$ or $M^\mathbb{Z}$) is to define entropy by  considering first a finite box of size, let's say $n$, and then take the limit on the size of the box: the thermodynamic limit. The probability on the finite box $M^n$ has no dynamical content. On the limit, when $n \to \infty$,  it may have  dynamical content (where the dynamics of shift corresponds to translation in the lattice $M^\mathbb{N}$ or $M^\mathbb{Z}$).  We say  in this case  that the entropy was obtained via finite partitions. In this setting, 
probabilities maximizing pressure are obtained in a similar way, like via the limit $\frac{e^{- H}\, dP}{\int e^{-H} \, d P}$, $n \to \infty$, where the Hamiltonian $H$ is in some way defined on each box of size $n$.
The procedure is different  in Thermodynamic Formalism, where you work   primarily with the Shannon--Kolmogorov entropy on the lattice $M^\mathbb{N}$ or $M^\mathbb{Z}$ (which has dynamical content) for getting shift invariant probabilities that maximize pressure. This entropy can be estimated by a version of the  Rokhlin Formula (see \cite{LMMS}). 
The Ruelle operator also played an important role in our definition of entropy. 
Both concepts are  linked
in a natural and fundamental way (see  \cite{LMMS}, or section 4 in \cite{BKL} for the classical thermodynamic formalism case).

In \cite{BKL} the authors show a relation of the entropy presented here
with Lyapunov exponents, and this is a clear indication of its dynamical nature.

Below we will present some clarifications on which directions our work is related to relevant issues in the area related to quantum entropy. 

First of all, is needed to say that the von Neumann entropy, which is given by  - trace $( \rho \log \rho)$,  in the same way as the expressions $- \sum_{i=1}^dp_i \log p_i$, or $\int \log f(x) f(x)   dx$, where $f$ is positive and
 $\int f(x) dx=1$, are not exactly  dynamical entropies (at least from our point of view).
 
 Quantum entropies with dynamical content were considered in a large number of papers and books for several decades.
 We believe our point of view does not coincide exactly (as far as we know) with the quite important results on the topic we describe next.

 In \cite{Araki} and \cite{Araki1} H. Haraki considers the relative
entropy which can be defined   for arbitrary normal states on a von Neumann algebra. As it is a relative  entropy is different from ours.

A very well know version is  the dynamical entropy of $C^*$-algebras and von Neumann
algebras of 
A. Connes, H. Narnhofer, and W. Thirring (see \cite{CNT}); as far we understand is in ''some sense  based'' on the principle of dynamic partitions.
 
 L. Accardi, A. Souissi and  E. Soueidy in \cite{Acar} consider a Quantum version of Markov Chains 
 which is in ''some sense''  based' on the principle of dynamic partitions. It is different from ours.
 
 R. Alicki and M. Fannes in  \cite{AF1}  considers the concept of  quantum dynamical entropy from different points of view: section 12 considers entropy production; section 13.1 consider the case of the quantum  cat map;  section 13.2 consider noncommutative 
 Lyapunov exponents and the Ruelle inequality (the dynamics are associated with the continuous-time  semigroup generated by the Laplacian in a compact Riemannian manifold); section  13.3 is devoted to quasi-free fermionic dynamics. All of them are different from ours.
 
 The setting of \cite{slomczynski2003dynamical} which considers iterated function systems and Markov operators  is the point of view closer to our work. But this reference does not consider the variational principle of pressure neither a version of the Ruelle operator. Results in  \cite{baraviera2010thermodynamic}, \cite{BLLT2} and \cite{BLLT1} addressed these topics and they  were generalized here.

The book \cite{Petz} consider the  relative {\it von Neumann entropy} in Quantum information with a view to some applications like the  Quantum Stein Lemma, Quantum Chernoff bounds, and Quantum Fisher information.

T. Sagawa in \cite{Sa} consider the relative entropy of von Neumann and questions related to the second law of Thermodynamics and majorization: what happens with the value of the entropy of a density 
matrix after the iteration by a quantum channel? The book \cite{Petz} addresses preliminarily   the question of majorization when a matrix is applied on a finite probability (\cite{LoRu} consider a similar problem considering the iteration of the dual of the Ruelle operator and not a matrix).  Maybe a future work could be to analyze majorization under  the context of the present paper.
 
  C. Pinzari, Y. Watatani, and K.  Yonetani in \cite{Pinzari} consider entropy and a variational principle of entropy  from the point of view of $C^*$-algebras. A version of the Perron-Frobenius theorem was used as an important tool for analyzing KMS states for some interesting examples arising from subshifts in symbolic dynamics. The relationship between the Voiculescu topological entropy 
and the topological entropy of the associated subshift is studied. 
In the case of the Cuntz-Krieger algebras, explicit construction of the state 
of maximal entropy was done. We understood that the  space of symbols (the alphabet) considered in \cite{Pinzari} is finite. Our results correspond to the case where  the alphabet (in some sense the support of the {\it a priori} probability $\mu$)
can be uncountable.

  In \cite{Pinzari1} the variational principle of pressure is considered by D. Kerr and C.  Pinzari. 
They introduce a notion of pressure for a selfadjoint element in a $C^*$-algebra, adapting Voiculescu's formulation of topological entropy for a nuclear $C^*$-algebra (see \cite{Voi} and \cite{Stor}).  The variational inequality holds for the Connes-Narnhofer-Thirring entropy.  They also introduce the concept of local state approximation entropy which is different from our definition  of entropy.


I. Nechita and C. Pellegrini addressed 
questions related to generic properties for quantum channels. In \cite{NePe} the authors show that for a fixed density matrix $\beta:\mathbb{C}^n \to \mathbb{C}^n$, the existence of a set of full measure for the Haar measure, on the set of unitary operator $U:\mathbb{C}^n \otimes \mathbb{C}^n \to \mathbb{C}^n \otimes \mathbb{C}^n $, satisfying the property that
for the associated quantum channel $Q \to \Phi(Q)  =Tr_{2} (U (q \otimes \beta) U^*)$ there exists a unique fixed point. In \cite{LoSe} the authors show that,  in fact, there exists an open and dense set of unitary operators $U$ with such property.

\medskip

A final remark: our main theorems   considered the case of the $C^*$-algebra of matrices $M_k$ and a natural question
is if our proofs can  be implemented for a general $C^*$-algebra?
Several results for completely positive maps that were used here are also known in a more general scope. This eventual extension would involve several issues that by their nature would be much more complex; in its generality would encompass - in a sense - the classical thermodynamic formalism
for potentials that depends on an infinite number of coordinates. The main eigenfunction for 
the Ruelle operator  of a continuous potential  may not exist; the existence  requires the use of the Holder regularity of
the potential. For the Markov case, the Perron Theorem provides similar
results without further hypotheses due to the fact that a potential that depends on two coordinates is automatic of Holder class. We leave the question related to the general $C^*$-algebra  for future work.

\end{document}